\documentclass[aap,MSNbibl,seceqn,dvips]{arximspdf}
\usepackage{graphicx}
\doi{10.1214/12-AAP893} %kopijuoti is PTS
\volume{23}
\issue{5}
\pubyear{2013}
\firstpage{2020}
\lastpage{2052}

\makeatletter

\newcommand{\rmP}{\mathrm{P}}

\newtheorem{theorem}{Theorem}[section]
\newtheorem{lemma}[theorem]{Lemma}
\newtheorem{corollary}[theorem]{Corollary}
\newtheorem{proposition}[theorem]{Proposition}

\newproclaim{example}[theorem]{Example}
\newproclaim{definition}[theorem]{Definition}
\newproclaim{assumption}[theorem]{Assumption}
\newproclaim{remark}[theorem]{Remark}

\newcommand{\R}{{\mathbb R}}
\newcommand{\F}{{\mathbb F}}
\newcommand{\Prob}{\mathbb{P}}
\newcommand{\E}{\mathbb{E}}
\newcommand{\sS}{{\mathcal S}}
\newcommand{\sL}{{\mathcal L}}
\newcommand{\armin}{\mathop{\arg\min}}
\newcommand{\armax}{\mathop{\arg\max}}
\newcommand{\linf}{\mathop{\lim\inf}}
\newcommand{\Fhat}{\hat{F}}

\makeatother

\begin{document}
\begin{frontmatter}

\title{Maximizing functionals of the maximum in the Skorokhod
embedding problem and an application to variance swaps}
\runtitle{Skorokhod embeddings and variance swaps}

\begin{aug}
\author[A]{\fnms{David} \snm{Hobson}\corref{}\ead[label=e1]{d.hobson@warwick.ac.uk}}
\and
\author[A]{\fnms{Martin} \snm{Klimmek}\thanksref{t1}\ead[label=e2]{m.klimmek@warwick.ac.uk}}
\runauthor{D. Hobson and M. Klimmek}
\affiliation{University of Warwick}
\address[A]{Statistics Department\\
University of Warwick\\
CV47AL, Coventry\\
United Kingdom\\
\printead{e1}\\
\hphantom{E-mail: }\printead*{e2}} %adresu isvedimo komanda gale!
\end{aug}

\thankstext{t1}{Supported from EPSRC funding.}

% HISTORY:
\received{\smonth{2} \syear{2012}}
\revised{\smonth{9} \syear{2012}}

% ABSTRACT
%
\begin{abstract}
The Az\'ema--Yor solution (resp., the Perkins solution) of the Skorokhod
embedding problem has the property that it maximizes (resp., minimizes)
the law of the maximum of the stopped process. We show that these
constructions have a wider property in that they also maximize (and
minimize) expected values for a more general class of bivariate
functions $F(W_\tau,S_\tau)$ depending on the joint law of the stopped
process and the maximum. Moreover, for monotonic functions $g$, they
also maximize and minimize $\E[\int_0^\tau g(S_t) \,dt ]$ amongst
embeddings of $\mu$, although, perhaps surprisingly, we show that for
increasing $g$ the Az\'ema--Yor embedding minimizes this quantity, and
the Perkins embedding maximizes it.\vspace*{1pt}

For $g(s)=s^{-2}$ we show how these results are useful in calculating model
independent bounds on the prices of variance swaps.

Along the way we also consider whether $\mu_n$ converges weakly to
$\mu$ is a sufficient condition for the associated Az\'ema--Yor and
Perkins stopping times to converge. In the case of the Az\'ema--Yor
embedding, if the potentials at zero also converge, then the stopping
times converge almost surely, but for the Perkins embedding this
need not be the case. However, under a further condition on the
convergence of atoms at zero, the Perkins stopping times converge in
probability (and hence converge almost surely down a subsequence).
\end{abstract}

% KEYWORDS

\end{frontmatter}

%s1 #&#
\section{Introduction}
\label{sec-intro}

Let $W = (W_t)_{t \geq0}$ be Brownian motion, null at 0, and $\mu$ a
centered probability measure. Then the Skorokhod embedding problem (SEP)
(Skorokhod~\cite{Skorokhod65}) is to find a stopping time $\tau$ such
that the stopped process satisfies $W_\tau\sim\mu$. There are many
classical solutions to this problem (for a survey,
see Ob{\l}{\'o}j~\cite{Obloj04}), and further solutions continue to
appear in
the literature, including most recently Hirsch et
al.~\cite{HirschProfetaRoynetteYor10}.
Further impetus to the investigation of old and new solutions is derived
from
the connections between solutions of the SEP and model independent
bounds for the prices of options; for
a survey, see Hobson~\cite{Hobson10}.

Given the multiplicity of solutions to the SEP, it is natural to search
for embeddings with additional optimality properties.
In particular, if $\Psi$ is a
functional of the stopped Brownian path $(W_t)_{0 \leq t \leq\tau}$,
then these constructions aim to maximize $\Psi$ over (a suitable
subclass of)
embeddings of $\mu$. For example, if $F$ is an increasing function, and
$S_t = \sup_{s \leq t} W_s$, then the Az\'ema--Yor
solution~\cite{AzemaYor79a}
maximizes $\E[F(S_\tau)]$ over uniformly integrable embeddings, and the
Perkins embedding~\cite{Perkins86} minimizes the same quantity.

Our goal in this paper is to extend this result to functions $F =
F(W_\tau,S_\tau)$. Then, subject to regularity conditions, our first
result (Theorem~\ref{thmmain1a}) is that:
\begin{quote}
Suppose
$F_s(w,s)/(s-w)$ is monotonic decreasing in $w$. Then
$\E[F(W_\tau,S_\tau)]$ is minimized (resp., maximized) over
uniformly integrable embeddings
$\tau$ of $\mu$ by the Az\'ema--Yor (resp., Perkins) embedding.
\end{quote}

This result is a tool in the derivation of our second result,
Theorem~\ref{thmmain2}, which, again subject to regularity
conditions is as follows:
\begin{quote}
Suppose $g$ is increasing. Then
$\E[\int_0^\tau g(S_u) \,du]$ is minimized (resp., maximized) over
uniformly integrable embeddings
$\tau$ of $\mu$ by the Az\'ema--Yor (resp., Perkins) embedding.
\end{quote}

One approach to finding extremal values of $\E[F(W_\tau,S_\tau)]$ is to
utilize the work of Kertz and R\"osler~\cite{KertzRosler90},
Vallois~\cite{Vallois94} and Rogers~\cite{Rogers93} who
characterize the
possible joint laws of $(W_\tau,S_\tau)$. These characterizations
take the
form of constraints on the possible laws of $(W_\tau,S_\tau)$, but
that still
leaves our problem as a constrained optimization problem. In fact,
there are
parallels between equation 3.2 of Theorem~3.1 of Rogers \cite
{Rogers93}, and
some of the quantities that arise in our study (see
Remark~\ref{ridearelated}), but we shall not make direct use of this
connection.

At first sight the second result above may appear counterintuitive.
After all, for
increasing $g$ the
Az\'ema--Yor embedding maximizes the law of $g(S_\tau)$ so one might also
expect it to maximize the law of $\int_0^\tau g(S_u) \,du$. However, the
exact opposite is true, and the Az\'ema--Yor embedding minimizes the expected
value
of this quantity. We return to this issue in
Remark~\ref{remcounterintuitive}, where we explain this phenomenon.

One of our tools for solving the above problems is to solve the problem
in the case where $\mu$ has bounded support, and to approach the
case of a general measure by approximation. In order to carry out this
program we need to analyze when and whether convergence of probability
measures is sufficient to guarantee that the associated Az\'ema--Yor and
Perkins embeddings converge. This proves to be a delicate question.
Under the additional (and necessary) hypothesis that $\int_{\R} |x|
\mu_n(dx) \rightarrow\int_{\R} |x| \mu(dx)$, then indeed the
Az\'ema--Yor embedding of $\mu_n$ converges almost surely to the
Az\'ema--Yor embedding of $\mu$. However, this need not be the case for
the Perkins embedding, and the sequence of Perkins embeddings of $\mu_n$
may fail to converge on an almost sure basis.

We note that although the focus in this paper is on functionals
involving the running maximum, there is a
parallel set of results for functionals involving the running minimum
process. The corresponding results
can be easily proved by following the proofs given for the maximum and
making the appropriate changes.
Alternatively, given a Brownian motion $W$ and a centered target law
$\mu$, let $\tilde{\mu}$ be the
measure $\mu$ reflected around zero. Then, with $I_t = \inf_{s \leq
t} W_s$, the problem of minimizing
$\E[F(W_\tau,I_\tau)]$ over embeddings $\tau$ of $\mu$ is
equivalent to minimizing
$\E[F(-{W}_{\tilde{\tau}},-{S}_{\tilde{\tau}})]$, over embeddings
$\tilde{\tau}$ of $\tilde{\mu}$. See the
next section and Section~\ref{ssecvswap} for calculations along these
lines.\vspace*{-2pt}

%s2 #&#
\section{A variance swap on squared returns}
\label{sec-squaredreturns}

The original motivation for our study came from financial mathematics
and the pricing of variance swaps, and one of the contributions of this article
is to establish a link between variance swap bounds and Skorokhod embedding
theory. The implications of this connection are the subject of related work
\cite{HobsonKlimmek10}. Informed by the results
presented here, but necessarily using different methods, Hobson and Klimmek
\cite{HobsonKlimmek10} show how to construct
model-independent bounds and hedging strategies for a general family of
variance swaps. In this section we outline the link between variance swaps
and the second result from the
\hyperref[sec-intro]{Introduction}.\looseness=1

Let $X=(X_t)_{0 \leq t \leq T}$ represent the discounted price of a
financial asset. Under the assumption of no-arbitrage, there exists a
measure under which $X$ is a (local)-martingale. We may suppose that
there exists a filtered probability space $(\Omega, {\mathcal F},
{\mathbb F}, \Prob)$ such that $B$ is a ${\mathbb F}$-Brownian and such
that $X_t=B_{A_t}$ for a (possibly discontinuous) time-change $t
\rightarrow A_t$, null at 0. (If $X$ is continuous, then the
existence of such a time-change is guaranteed by the
Dambis--Dubins--Schwarz theorem, and in general the existence is
guaranteed by Monroe~\cite{Monroe78}, Theorem~2.) Since $X$ is a
nonnegative price process we suppose it has starting value $X_0= B_0
=x_0>0$.

Now suppose that we know the prices of put and call options with
maturity~$T$. Knowledge of put and call option prices with expiry time
$T$ is equivalent to knowledge of the marginal law of process at time
$T$; see Breeden and Litzenberger~\cite{BreedenLitzenberger78}.
Suppose that $X_T \sim\mu$ and that $\mu$ is centered at $x_0$, and has
support in $\R^+$.
%(then $X$ is a true martingale and not a strict local martingale).
We will determine bounds for the fair value of a variance
swap given the terminal law $\mu$. Note that if $X_T \sim\mu$, then
$A_T$ is a solution of the Skorokhod embedding problem for $\mu$ in $B$.

Following Demeterfi et al.~\cite{DemeterfiDerman99} we define the
pay-out $V= V((X_s)_{0 \leq s \leq T})$ of an idealized
variance swap as
%
%e2.1 #&#
\begin{equation}
\label{eqidealised} V_T =\int_0^T
\frac{d[X,X]_t}{(X_{t-})^2}=\int_0^T \biggl(
\frac{dX^c_t}{X_{t-}} \biggr)^2+\sum_{0 \leq t \leq T}
\biggl(\frac{\Delta X_t}{X_{t-}} \biggr)^2,
\end{equation}
where $\Delta X_t = X_t - X_{t-}$, and $X^c$ is the continuous part of
$X$.

Let\vspace*{1pt} $A^c$ be the continuous part of $A$. Note that
$dA^c_t=(dX_t^c)^2=d[X,X]^c_t$.
Let $S^X=(S^X_t)_{t \geq0}$ (resp., $S^B$)\vspace*{1pt} be the process of
the running maximum of $X$ (resp.,~$B$), and let $I^X$
(resp., $I^B$)
denote
the corresponding infimum.
Then we have
%$I_{A_t} \leq I^X_t \leq
$X_t \leq S^X_t \leq S^B_{A_t}$, and it follows that path-by-path
with $\Delta B_{A_t} = B_{A_t} - B_{A_{t-}}$ that
%
%e2.2 #&#
\begin{eqnarray}
\label{eqchap3contmon1} V_T &\geq& \int_0^T
\frac{d[X,X]^c_t}{(S^X_{t-})^2} + \sum_{0 \leq t \leq T} \biggl(
\frac{\Delta X_t}{S^X_{t-}} \biggr)^2
\nonumber\\[-8pt]\\[-8pt]
&\geq& \int_0^T \frac{d A_t^c}{(S^B_{A_{t-}})^2} + \sum
_{0 \leq t
\leq T} \biggl(\frac{\Delta B_{A_t}}{S^B_{A_{t-}}} \biggr)^2.
\nonumber
\end{eqnarray}
We suppose that $X$ has a second moment. Then
$(X_{t})_{0 \leq t \leq T}$ is a square-integrable martingale and we
find that
%
%e2.3 #&#
\begin{eqnarray}
\label{eqchap3contmon2} \E\biggl[\int_0^T
\frac{dA_t^c}{(S^B_{A_{t-}})^2}+\sum_{0 \leq t
\leq
T} \biggl(\frac{\Delta B_{A_t}}{S^B_{A_{t-}}}
\biggr)^2 \biggr] &=& \E\biggl[\int_0^T
\frac{dA^c_{t}+\Delta A_t}{(S^B_{A_{t-}})^2} \biggr]
\nonumber\\
&=& \E\biggl[\int_0^T \frac{dA_t}{(S^B_{A_{t-}})^2}
\biggr]
\\
&\geq& \E\biggl[\int_0^{A_T} \frac{du}{(S^B_u)^2}
\biggr].
\nonumber
\end{eqnarray}

We say that $\tau$ is an embedding of $\mu$ if $\tau$ is a stopping
time for
which $B_\tau$ has law~$\mu$ [we write $B_\tau\sim\mu$ or $\mu=
\sL(B_\tau)$]. Let $\sS\equiv\sS(B,\mu)$ be the set of stopping
times which
embed $\mu$, and let $S_{\mathrm{UI}} = \sS_{\mathrm{UI}}(B,\mu)$ be the subset of
$\sS(B,\mu)$ for
which $(B_{t \wedge\tau})_{t \geq0}$ is uniformly integrable.
The inequalities above imply that the fair value of $V_T$ is
bounded below by
%
%e2.4 #&#
\begin{equation}
\label{eqvarlower} \inf_{\tau\in\sS_{\mathrm{UI}}(B,\mu)} \E\biggl[\int_0^\tau
\frac
{du}{(S^B_u)^2} \biggr].
\end{equation}
Similarly, using the inequality $I^B_{A_t} \leq I^X_t \leq X_t$ we find
that the
fair value of $V_T$ is bounded above by
%
%e2.5 #&#
\begin{equation}
\label{eqvarupper} \sup_{\tau\in\sS_{\mathrm{UI}}(B,\mu)} \E\biggl[\int_0^\tau
\frac
{du}{(I^B_u)^2} \biggr].
\end{equation}
This problem can be converted into a problem concerning the maximum
$S^B$ by a
reflection argument; see Section~\ref{ssecvswap}.

Now let $G(b,s)=\frac{(s-b)^2}{s^2}$. Then by It\^{o}'s lemma,
\[
G\bigl(B_\tau,S^B_\tau\bigr)=G(0,0)+\int
_0^\tau\frac{du}{(S^B_u)^2} - \int
_0^\tau\frac{2(S^B_u-B_u)}{(S^B_u)^{2}} \,dB_u.
\]
It follows that
if $\int_0^{\tau\wedge t} {2(S^B_u-B_u)}{(S^B_u)^{-2}} \,dB_u$ is a uniformly
integrable
martingale, then
\[
\E\biggl[\int_0^\tau\frac{du}{(S^B_u)^2}
\biggr]=\E\biggl[\frac{(S^B_\tau-B_\tau
)^2}{(S^B_\tau)^2} \biggr],
\]
and the question of bounding the fair value of $V_T$ is transformed
into a
question of maximizing or minimizing expressions of the form
$\E[F(B_\tau,S_\tau)]$ over embeddings of $\mu$. We
return to the calculation of the variance swap bound in Section \ref
{ssecvswap}.

%s3 #&#
\section{Preliminaries}
\label{sec-prelims}

Let $( \Omega, {\mathcal F}, {\mathbb F}, \Prob)$ be a filtered
probability space satisfying the usual conditions and supporting a
Brownian motion $W=(W_t)_{t \geq
0}$ with
\mbox{$W_0=0$}, and sufficiently rich that ${\mathcal F}_0$ contains a further
uniform random variable which is independent of $W$. Let $\mu$ be a
centered probability measure. To exclude trivialities we assume that
$\mu$ is not $\delta_0$, the unit mass at 0. We say that
$\tau$ is an embedding of $\mu$
if $\tau$ is a stopping time for which $W_\tau$ has law $\mu$ [we
write $W_\tau\sim\mu$ or $\mu= \sL(W_\tau)$] and we say
that $\tau$ is uniformly integrable if the family $(W_{t \wedge
\tau})_{t
\geq0}$ is uniformly integrable.

Let $\sS\equiv\sS(W,\mu)$ be the set of stopping times which embed
$\mu$, and let
$\sS_{\mathrm{UI}} \equiv\sS_{\mathrm{UI}}(W,\mu)$ be the subset of $\sS(W,\mu)$
consisting of uniformly
integrable stopping times.
For $\sS_{\mathrm{UI}}(W,\mu)$ to be nonempty we must have that $\mu$ is centered
[i.e., $\int_{\R} |x| \mu(dx) < \infty$ and $\int_{\R} x \mu(dx) =0$].
In this context (Brownian motion and centered
target laws) a result of Monroe~\cite{Monroe72} gives that a stopping
time is uniformly integrable if and only if it is minimal (in the sense
that if $\tau$ is minimal and $\sigma\leq\tau$ with $W_\sigma\sim
W_\tau$, then $\sigma\equiv\tau$ almost surely). The class of minimal
stopping
times is a natural class of ``good'' (in the sense of small) stopping
times.

For the Brownian motion $W$, started at 0, we write $H_{x}$ for the
first hitting time of $x$, and for a set $A$, $H_A = \inf\{ u \geq0\dvtx
W_u \in A \}$.

For a process $(Y_t)_{t \geq0}$ and a stopping time $\sigma$, we write
$Y^\sigma= (Y^\sigma_t)_{t \geq0}$ for the stopped process
$Y^\sigma_t = Y_{\sigma\wedge t}$.

Given a centered probability measure $\mu$, let $X_\mu$ be a random
variable
with law~$\mu$, and define $C(x) \equiv C_\mu(x) = \E[(X_{\mu} - x)^+]$
and $P(x) \equiv P_\mu(x) = \E[(x - X_\mu)^+]$. Then $C$ and $P$ are
monotonic convex functions with $C(0) = P(0)$. Then $U(x) = U_\mu(x) =
\E[|X_{\mu} - x|] = C(x)+P(x)$ is (minus) the potential associated with
$\mu$. Conversely any convex function $U$ with $\lim_{x \rightarrow
\pm\infty} (U(x)- |x|) = 0$ is the potential of some centered
probability measure $\mu$ (Chacon~\cite{Chacon77}).\vadjust{\goodbreak}

If $\mu$ has an atom at zero, then we write $\mu^*$ for the measure
obtained by omitting the atom at 0, and then rescaling to get a
probability measure. Thus $\mu^*(A) = \mu(A \setminus\{0\})/(1-
\mu(\{0\}))$. Finally, we write $\hat{x} = \hat{x}_\mu$ for the
upper limit
on the support of $\mu$ [so $\hat{x}_\mu= \sup\{ x\dvtx  C_\mu(x)>0
\}$]
and
$\check{x} = \check{x}_\mu$ for the corresponding lower limit
$\check{x}_\mu= \inf\{ x\dvtx  P_\mu(x)>0 \}$.

%f1 #&#
\begin{figure}

\includegraphics{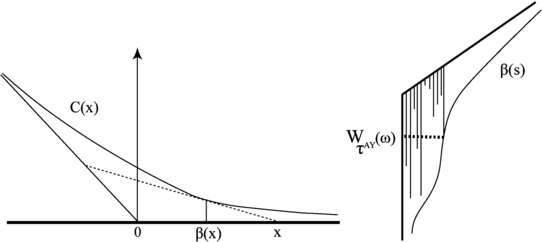}

\caption{For each $x$, the value of $\beta(x)$ is determined by finding
the tangent line to $C_\mu$ originating at $x$: $\beta(x)$ is the
horizontal co-ordinate of the point of contact between the tangent line
and $C_\mu$. [If $C_\mu$ includes a straight line section, then this
point of contact may not be uniquely defined, in which case we take
$\beta(x)$ to be the largest value of the horizontal co-ordinate at
which contact occurs.] The stopping time $\tau_\beta$ associated to this
construction is given by the first time that an excursion from the
maximum crosses below $\beta$.}\label{fig1}
\end{figure}

%s3.1 #&#
\subsection{The Az\'{e}ma--Yor solution}
\label{sec-azemayor}
For $x \geq0$, up to the upper limit on the support of $\mu$,
define $\beta= \beta_\mu$ by
%
%e3.1 #&#
\begin{equation}
\label{eqnbetadef} \beta(x) = \armin_{y<x} \frac{C_\mu(y)}{x-y}.
\end{equation}
Then $\beta$ is an increasing function with $\beta(x)<x$, see Figure
\ref{fig1}. Where the $\armin$ is not uniquely defined it is not
important which value we choose. However, we fix one by insisting that
$\beta$ is right-continuous, or equivalently by choosing the largest
value for which the minimum is attained. Observe that at $x=0$, $\beta$
takes the value of the infimum of the support of $\mu$. For $x$ equal
to, or to the right of, the upper limit on the support of $\mu$ we set
$\beta(x)=x$.

For an increasing function $\beta\dvtx \R^+ \mapsto\R$ with $\beta(x)
\leq
x$ let $\tau_\beta$ be given by
%
%e3.2 #&#
\begin{equation}
\label{eqntaubetadef} \tau_\beta= \inf\bigl\{u\dvtx  W_u
\leq\beta(S_u) \bigr\}.
\end{equation}
Then $\tau^{\mathrm{AY}} \equiv\tau^{\mathrm{AY}}_\mu$,
the Az\'ema--Yor stopping time for $\mu$, is given by
$\tau^{\mathrm{AY}}_\mu\equiv\tau_{\beta_\mu}$. Thus
we have $\tau_{\beta_\mu} \in\sS_{\mathrm{UI}}(W,\mu)$, and moreover, for $F$
increasing,
$\tau_{\beta_\mu}$ maximizes $\E[F(S_\tau)]$ over $\tau\in
\sS_{\mathrm{UI}}(W,\mu)$
(Az\'ema--Yor~\cite{AzemaYor79a,AzemaYor79b}, Rogers~\cite{Rogers89}).

Note that $\tau_{\beta_\mu}$ does not maximize this quantity over
\textit{all}
embeddings, but it does give the maximum over uniformly integrable
(i.e., minimal) embeddings.

Let $b\equiv b_\mu$ be the right-continuous inverse to $\beta$. Then $b$
is the barycenter function and for $x < \hat{x}_\mu$, $b(x)$ is given by
%
%e3.3 #&#
\begin{equation}
\label{barycentre} b(x) = \E[X_\mu| X_\mu\geq x].
\end{equation}
The barycentre $b(x)$ is defined
up to the
upper limit of the support
of $\mu$ and
is a nonnegative, nondecreasing function with $b(x) \geq x$.
We set $b(x)=x$ for $x\geq\hat{x}_\mu$. (The reverse barycentre
$\check{b}(x) = \E[X|X \leq x]$ is defined analogously to the
barycentre.)

It is more standard to define the barycenter function as in
(\ref{barycentre}) and to set $\beta$ to be
the inverse barycenter function, but the two approaches are equivalent,
and our approach via potentials allows for a unified treatment with the
Perkins
construction in the next section.

If $\mu$ has an interval with no mass, then $b$ is constant over that
interval, and $\beta$ has a jump. If $\mu$ has an atom at $x$ then $b$
has a jump at $x$ [unless the atom is at the upper limit $\hat{x}$ of
the support of $\mu$ in which case $b(\hat{x})=\hat{x}$] and $\beta
$ is
constant over a range of $s$. From the definition of $\tau_\beta$ [see
(\ref{eqntaubetadef})] and
excursion theory (see Rogers~\cite{Rogers93}, equation 2.13), we
have
%
%e3.4 #&#
\begin{equation}
\label{eqrogersfact} \exp\biggl( - \int_0^s
\frac{dr}{r - \beta(r)} \biggr) = \Prob(S_{\tau_\beta} \geq s)
\end{equation}
and then also $\Prob(S_{\tau_\beta} \geq s)
= \Prob(W_{\tau_\beta} \geq\beta(s)) = \mu(\beta(s),\infty)$.
Note that it does not matter which convention we use for $\beta(s)$
here since
$\mu$ places no mass on $(\beta(s-), \beta(s+))$.
%
%ex3.1 #&#
\begin{example}
\label{excalcuniform} If $\mu=U[-1,1]$, then $C_{\mu}(x) = (x-1)^2/4$
and $P_{\mu}(x) = (x+1)^2/4$
(at least for $-1 = \check{x}_\mu\leq x \leq\hat{x} = 1$).
Then
the barycenter
function is given by $b(x)=({x+1})/{2}$ for $-1 \leq x \leq1$ and
hence
$\beta(s)=2s-1$ for $0 \leq s \leq1$. It follows that
$S_{\tau^{\mathrm{AY}}_\mu} \equiv b (W_{\tau^{\mathrm{AY}}_\mu})$ is uniformly
distributed on $[0,1]$.
\end{example}

%le3.2 #&#
\begin{lemma}
\label{lemintegrabilitybeta}
If $\mu$ places mass on $(x,\infty)$, then
$(r - \beta(r))^{-1}$ is integrable over $[0,x]$.
\end{lemma}
\begin{pf}
This follows immediately from (\ref{eqrogersfact}) and
$\Prob(S_{\tau_\beta}
\geq x)
\geq\Prob(W_{\tau_\beta} \geq x) > 0$.
\end{pf}

%s3.2 #&#
\subsection{The Perkins solution}
\label{secperkins}

\begin{figure}

\includegraphics{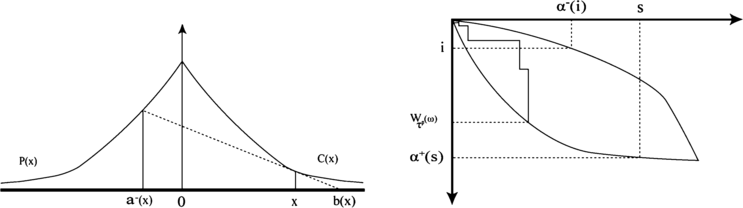}

\caption{Suppose that $\mu$ has no atoms. Then for $x>0$,
$a^-(x)$ is the horizontal co-ordinate of the point where the tangent
line to $C$ at $(x,C(x))$ intersects
with $P$. Alternatively, it is the horizontal co-ordinate of the point
where the tangent line to $C$
emanating from $(b(x),0)$ intersects with $P$. [We may instead consider
the inverse $\alpha^-$ of
$a^-$:  for $y<0$, $\alpha^-(y)>0$ is the horizontal co-ordinate of
the point such that the tangent to $C$ at $\alpha^-(y)$ crosses $P$ at
$(y,P(y))$.]
These definitions extend naturally to the case where the convex
function $C$ has
kinks or straight-line segments.
Similarly, $a^+(z)$ is found by drawing tangents to $P$ emanating
from the reverse barycenter function evaluated at $z<0$
and determining intersection points with $C$.
The stopping rule associated with this construction is to stop the
Brownian motion
when its running maximum or minimum exit the region determined by
$\alpha^+$ and $\alpha^-$.}\label{fig2}
\end{figure}

For $x>0$ define $\alpha^+_\mu= \alpha^+\dvtx  \R_+ \rightarrow\R_-$ by
%
%e3.5 #&#
\begin{equation}
\label{alpha+} \alpha^+(x) = \armin_{y<0} \frac{C_\mu(x)-
P_{\mu}(y)}{x-y}
\end{equation}
and for $x<0$ define $\alpha^-_\mu= \alpha^-\dvtx  \R_- \rightarrow\R_+$ by
%
%e3.6 #&#
\begin{equation}
\label{alpha-} \alpha^-(x) = \armax_{y>0} \frac{P_\mu(x)-
C_{\mu}(y)}{y-x}.
\end{equation}
Then $\alpha^{\pm}$ are monotonic functions, see Figure~\ref{fig2}. If
the $\armin$ (resp., the $\armax$) is not uniquely defined, we take the
largest value (in modulus) for which the minimum (resp., the maximum)
is attained; in this way $\alpha^{+}\dvtx  \R_{+} \mapsto\R_{-}$ is
right-continuous and $\alpha^-$ is left-continuous. Again, none of the
subsequent analysis will depend on this convention. For convenience we
will sometimes write $\alpha$ as shorthand for $\alpha^{\pm}$.

%f2 #&#

If $P_{\mu}$ (resp., $C_\mu)$
is differentiable at $\alpha^+(x)$ [resp., $\alpha^-(x)$], then
$\alpha^+(x)$
[resp., $\alpha^-(x)$] satisfies
%
%e3.7 #&#
\begin{equation}
\label{eqalphaplussolves} \frac{C_\mu(x)- P_{\mu}(\alpha
^+(x))}{x-\alpha^+(x)} = P'_{\mu
}
\bigl(\alpha^+(x)\bigr)
\end{equation}
[resp., $P_\mu(x)- C_{\mu}(\alpha^-(x)) =
C'_{\mu}(\alpha^-(x))(x-\alpha^-(x))$].

Let $a^\pm$ be the inverse to $\alpha^{\pm}$ and let $\bar{a}(w) =
w$ for
$w>0$ and $\bar{a}(w)=a^+(w)$ for $w<0$.
Recall the definition of $I$ as the infimum process for $W$ so that
$I_t = \inf_{s \leq t} W_s$.

For a pair of monotonic functions $\alpha^{+}\dvtx  \R_{+} \mapsto
\R_{-}$ (nonincreasing) and $\alpha^{-}\dvtx  \R_{-} \mapsto
\R_{+}$ (nondecreasing) define the stopping time
\[
\tau_\alpha= \inf\bigl\{u>0\dvtx  W_u \leq
\alpha^+(S_u) \bigr\} \wedge\inf\bigl\{u>0\dvtx  W_u \geq
\alpha^-(I_u) \bigr\}.
\]

Suppose $\mu$ does not have an atom at zero.
Then the Perkins~\cite{Perkins86} embedding $\tau^{\rmP } \equiv
\tau^{\rmP }_{\mu} \equiv\tau^{\rmP }(\mu)$ is given by $\tau^{\rmP }_\mu=
\tau_{\alpha_{\mu}}$.

If $\mu$ has an atom at zero,
%then define $\mu^*$ via
%$\mu^*(A) = \mu(A \setminus\{0\})/(1 - \mu(\{0\})$. Then $\mu^*$ is
%$\mu$
%with the mass at zero removed, rescaled to be a probability measure.
then we use independent randomization to
set $\tau^{\rmP } = 0$ with probability $\mu(\{0 \})$;
and otherwise $\tau^{\rmP } = \tau_{\alpha_{\mu}}$.
More precisely, in the case where $\mu$ has an atom at zero
we set the Perkins embedding to be
\[
\tau^{\rmP } = \cases{0, &\quad if $Z\leq\mu({0})$,
\cr
\tau_{\alpha_\mu}, &\quad
if $Z> \mu({0})$,}
\]
where $Z$ is a uniform random variable which
is measurable with respect to
${\mathcal F}_0$.
Here
$\alpha^{\pm}_\mu$ are the quantities defined in
(\ref{alpha+}) and (\ref{alpha-}) for $\mu$. Note that if $\mu^*$ is
obtained from
$\mu$ by removing any mass at zero, and rescaling to give a probability
measure, then although $C_{\mu^*}$ and $P_{\mu^*}$ are
scalar multiples of $C_\mu$ and $P_\mu$, respectively, nonetheless we
have $\alpha^{\pm}_{\mu^*} \equiv\alpha^{\pm}_{\mu}$.

Note that if $\mu$ has an atom at zero, then we need ${\mathcal F}_0$
to be
nontrivial in order to be able to define the Perkins embedding. Note further
that since there are potentially many uniform random variables $Z$
which are
measurable with respect to ${\mathcal F}_0$, if $\mu(\{0\})>0$, then the
Perkins embedding is not unique. Sometimes it is convenient to
think about the Perkins embedding
associated with an identified ${\mathcal F}_0$ random variable $Z$, in which
case we write $\tau^{\rmP ,Z}_{\mu}$ instead of just $\tau^{\rmP }_{\mu}$.

The results of Perkins~\cite{Perkins86} show that $\tau^{\rmP }_\mu\in
\sS_{\mathrm{UI}}(W,\mu)$ and moreover, for $F$
increasing,
$\tau^{\rmP }$ minimizes $\E[F(S_\tau)]$ over $\tau\in\sS(W,\mu)$, and
not just $\sS_{\mathrm{UI}}(W,\mu)$ (Perkins~\cite{Perkins86}, although
the representation via (\ref{alpha+}) and (\ref{alpha-}) is due to
Hobson and Pedersen~\cite{HobsonPedersen02}).
%
%ex3.3 #&#
\begin{example} \label{excalcuniform2}
If $\mu=U[-1,1]$, then $P=P_{\mu}$ and $C=C_\mu$ are as given in
Example~\ref{excalcuniform} and, from (\ref{eqalphaplussolves}),
$\alpha^+(s)$ is the unique root of the equation
$P'(\alpha)(s-\alpha)=C(s)-P(\alpha)$.
It is easily verified that this root is given by
$\alpha^+(s) = s- 2\sqrt{s}$. Similarly, $\alpha^-(i) = i +
2 \sqrt{|i|}$. It can be shown that $\Prob(S_{\tau_\alpha} \geq s) =
\Prob(W_{\tau_\alpha} \geq s) = \Prob(W_{\tau_\alpha} \leq s - 2
\sqrt{s}) = 1 - \sqrt{s}$.
\end{example}
%
%ex3.4 #&#
\begin{example} \label{excalcpareto}
Notwithstanding the above example, in general it is difficult to derive
an explicit
form for the
stopping boundary associated with the Perkins stopping
time. Here we give a second example where analytic expressions,
albeit complicated ones, can be derived.

Suppose the target law is a
centered Pareto distribution with support $[-1,\infty)$ and density
function $f(x)=2(x+2)^{-3}$.
Then for $k \geq-1$,
$C(k) = (2+k)^{-1}$ and $P(k) = k + ({2+k})^{-1}$, and for $k<-1$,
$C(k)=-k$, $P(k)=0$.

Then, for the Az\'{e}ma--Yor embedding, $\beta$ solves $C(\beta) =
(s-\beta)
|C'(\beta)|$ and $\beta(s)=(s/2)-1$.

For the Perkins embedding, $\alpha^+(s)$ solves
$P'(\alpha^+)=({C(s)-P(\alpha^+)})/\break({s-\alpha^+})$, and we have (after
some algebra)
\[
\alpha^+(s) = \frac{-2s^2-5s+\sqrt{s^4+6s^3+12s^2+8s}}{2s-1+s^2}.
\]
The expression for $\alpha^-$ is $\alpha^-(i)=\frac{-3i-2i^2+\sqrt
{-(i^4+6i^3+12i^2+8i)}}{2i+1+i^2}$.\vadjust{\goodbreak}
\end{example}
If $\mu$\vspace*{1pt} has an interval in $\R_+$ (resp., $\R_-$) with no mass,
then $\alpha^-$ (resp., $\alpha_+$) has a jump (unless that
interval is contiguous with zero, in which case $\alpha^{\pm}$ starts at
a nonzero value).
If $\mu$ has an atom in $(0,\infty)$ [resp.,
$(-\infty,0)$], then $\alpha^-$ (resp., $\alpha_+$) is constant
over a range of values.
%
%le3.5 #&#
\begin{lemma}
\label{lemintegrabilityalpha}
Suppose $x>0$. If $\mu$ places mass on $[x,\infty)$, then
$(r - \alpha^+(r))^{-1}$ is integrable over $(0,x)$.
\end{lemma}
\begin{pf}
We have $(W_u \geq\alpha^+(S_u); \forall u \leq H_x) \supseteq
(\tau_\alpha\geq H_x) \supseteq(W_{\tau_\alpha} \geq x)$, and then
by excursion theory [recall (\ref{eqrogersfact})],
\[
\exp\biggl( - \int_0^x \frac{dr}{r - \alpha^+(r)}
\biggr) = \Prob\bigl(W_u \geq\alpha^+(S_u); \forall u
\leq H_x\bigr) \geq\mu\bigl([x,\infty)\bigr)>0.
\]
\upqed\end{pf}

%s4 #&#
\section{Convergence of measures and convergence of embeddings}
\label{secapp}

Let $(\mu_n)_{n \geq1}$ be a sequence of measures, and write $U_n$,
$\beta_n$ and $\alpha_n$ as shorthand for $U_{\mu_n}$, $\beta_{\mu_n}$
and~$\alpha_{\mu_n}$, with a similar convention for other functionals.

Suppose that, for each $n$, $\mu_n$ is centered
and that $(\mu_n)_{n \geq1}$
converges weakly to $\mu$, where $\mu$ is also centered. Then it does not
follow that
$U_{n} \rightarrow U_\mu$, nor that $\beta_n \rightarrow
\beta_\mu$, nor that $\alpha_n \rightarrow\alpha_\mu$. However, with
the
correct additional hypotheses, then these types of convergence are
equivalent.
%unnecessary now
%We begin with a result from convex analysis.
%Let $f\dvtx \R\rightarrow\R$
%be a convex function, let $f^c$ denote the convex conjugate of $f$, so
%that $f^c(y) = \sup_{x}
%$f$ at $x$ (defined as the set of points where the supremum above is
%attained). Note that $\partial f(x)$ is nonempty. In general $
%an interval and we write
%$\theta_f(x) = \inf\{ y \in\partial f(x) \}$ and
%$\psi_f(x) = \sup\{ y \in\partial f(x) \}$.
%Then $f$ is differentiable at $x$ if and only $\partial f(x)$ is single
%valued if and only if $\psi_f(x) = \theta_f(x)$.
%
%The
%following lemma is adapted from Rockafellar \cite Theorem
%24.5 {Rockafellar\dvtx 70}.
%
%Let $(g_n)$ be a sequence of
%finite convex functions on $\R$ converging pointwise to the convex
%function $g$.
%Let $(z_n)$ be a sequence of points converging to $z$.
%Given any $\varepsilon>0$ there exists an $N \in\N$ such that for $n \geq
%N$,
%(\theta_g(z)-\varepsilon,\psi_g(z)+\varepsilon).\]

Our first key result is the following.
%
%pr4.1 #&#
\begin{proposition} \label{pconvergence}
Let $(\mu_n)$ be a sequence of measures such that $\mu_n \Rightarrow
\mu$ and
$\E[|X_{\mu_n}|] \rightarrow\E[|X_\mu|]$. Then $b_n(x) \rightarrow
b(x)$ at continuity points $x < \hat{x}$ of $b$.
\end{proposition}
\begin{pf}
Chacon~\cite{Chacon77} shows that if $\mu_n \Rightarrow\mu$ and
$U_n(0) \rightarrow U(0)$, then $U_n \rightarrow U$ pointwise.
Since $C_n(x)=(U_n(x)+x)/2$ it follows trivially that $C_n \rightarrow
C$ pointwise, where
$C_n(x)=C_{\mu_n}(x)$ and $C(x)=C_{\mu}(x)$.

Recall that $x$ is a discontinuity point of $b$ if and only if there is an
atom of
$\mu$ at~$x$.
Suppose $x<\hat{x}$ is a continuity point of $b$. Then
(\ref{barycentre}) gives $b(x)=x+\frac{C(x)}{\mu([x,\infty))}$ and
\[
b_n(x)=x+\frac{C_n(x)}{\mu_n([x,\infty))} \rightarrow x+\frac
{C(x)}{\mu([x,\infty))} =
b(x).
\]
\upqed\end{pf}

\begin{corollary}
\label{corbtobeta}
Let $(\mu_n)$ be a sequence of measures such that $\mu_n \Rightarrow
\mu$
and
$\E[|X_{\mu_n}|] \rightarrow\E[|X_\mu|]$. Then $\beta_n(s)
\rightarrow
\beta(s)$ at continuity points $s < \hat{x}$ of $\beta$.
Moreover, if $\hat{x}<\infty$, then for
each $z> \hat{x}$,
$\lim\inf\beta_n(z) \geq\hat{x}$.
\end{corollary}
\begin{pf}
Since $b_n(\hat{x} - \varepsilon) < \hat{x}+\varepsilon$ for sufficiently
large
$n$ we have for these same $n$ that $\beta_{n}(\hat{x}+\varepsilon)
\geq
\hat{x} - \varepsilon$.
\end{pf}
%
%co4.3 #&#
\begin{corollary} \label{cazemaconv} Under the assumptions of
Proposition
\ref{pconvergence}, $\tau_{\beta_n} \rightarrow\tau_\beta$ almost
surely.\vadjust{\goodbreak}
\end{corollary}
\begin{pf}
Let $D$ be the set of discontinuity points of $\beta$.
If $ S_{\tau_\beta} \notin D$, then $W_{\tau_\beta} =
\beta(S_{\tau_\beta})$, and it follows that
\[
(\omega\dvtx  \tau_{\beta_n} \not\rightarrow\tau_\beta) \subseteq(
\omega\dvtx  S_{\tau_\beta} \in D ) \cup\bigl(\omega\dvtx  S_{\tau_\beta} \notin D,
W_{\tau_\beta} = \beta(S_{\tau_\beta}), \tau_{\beta_n} \not\rightarrow
\tau_\beta\bigr).
\]
For any stopping time $\sigma$ write:
let $H^\sigma_x = \inf\{ u \geq\sigma\dvtx  W_u = x \}$.

\textit{Case} 1: $\hat{x}=\infty$.
Note that since $\beta$ is increasing, $D$ is countable and
$\Prob(S_{\tau_\beta}
\in D )=0$.

First we argue that on $(\omega\dvtx  S_{\tau_\beta}=x)$ we have that for
sufficiently large $n$,
\mbox{$S_{\tau_{\beta_n}} \geq x$}: since there are only countably many
values of
$s<x$ on which the value of $W_u$ gets below $S_u=s$, and on each of
these excursions $W$ stays above $\beta(S)$, for sufficiently large $n$,
$W$ must stay above $\beta_n(S)$ also.

Hence $\liminf_n S_{\tau_{\beta_n}} \geq S_{\tau_\beta}$ almost
surely. Then on
$\{ \omega\dvtx  S_{\tau_\beta}= x \notin D, W_{\tau_\beta}=\beta(x) \}
$, we
have $ \tau_{\beta_n}(\omega) \rightarrow\tau_\beta(\omega)$ unless
$\inf\{ W_u\dvtx  \tau_{\beta} \leq u \leq H^{\tau_\beta}_{S_{\tau
_\beta}}
\} = W_{\tau_\beta} = \beta(x)$ and $\beta_n(x) < \beta(x)$. But, almost
surely, on any interval of positive length Brownian motion goes below
its starting value.
In particular, the set
$(\omega\dvtx  S_{\tau_\beta} \notin D, W_{\tau_\beta}=\beta(S_{\tau
_\beta}),
\tau_{\beta_n}
\not\rightarrow\tau_\beta)$ has probability zero.

\textit{Case} 2: \textit{$\hat{x}<\infty$ and $\mu( \{\hat{x} \}
)=0$.}
The only paths for which issues of convergence might be different to
the previous case are those for which $S_{\tau_\beta}=\hat{x}$. But
since $\mu$ has no atom at $\hat{x}$, $\Prob(S_{\tau_\beta} =
\hat{x}) =
\Prob(W_{\tau_\beta} = \hat{x})=0$ and $\tau_{\beta_n}
\rightarrow
\tau_\beta$ almost surely.

\textit{Case} 3: \textit{$\hat{x}<\infty$ and $\mu( \{\hat{x} \}
)>0$.}
In this case $\beta(\hat{x}-):= \lim_{y \uparrow\hat{x}} \beta
(y) <
\beta(\hat{x}) = \hat{x}$. We show that on the set
$(S_{\tau_\beta}=\hat{x})$
we
have $\lim\tau_{\beta_n} = \tau_\beta$, almost surely. Off the set
$(S_{\tau_\beta}= \hat{x})$
convergence follows exactly as in the previous cases.

First we argue that $\lim\sup_n S_{\tau_{\beta_n}} \leq\hat{x}$ almost
surely. Fix $z>\hat{x}$, then given $0<\varepsilon<z-\hat{x}$, there exists
$N$ such
that for $n \geq N$, $\beta_n(\hat{x}+\varepsilon)>\hat{x}-\varepsilon$.
Hence,
for sufficiently large $n$,
\[
\bigl( \omega\dvtx  S_{\tau_{\beta_n}}(\omega) \geq z \bigr) \subseteq\bigl
( \omega\dvtx
\inf\{W_u\dvtx  H_{\hat{x}+\varepsilon} \leq u \leq H_z \} \geq
\hat{x}-\varepsilon\bigr).
\]
But
\[
\Prob\bigl( \inf\{W_u\dvtx  H_{\hat{x}+\varepsilon} \leq u \leq
H_z \} \geq\hat{x}-\varepsilon\bigr) \leq\exp\biggl( - \int
_{\hat{x}+\varepsilon}^z \,\frac{dy}{y - (\hat{x}-\varepsilon)} \biggr) =
\frac{2 \varepsilon}{z -
\hat{x}+\varepsilon}.
\]
By choosing $\varepsilon$ small compared with $(z-\hat{x})$ we deduce that
$\lim\sup_n S_{\tau_{\beta_n}} \leq z$ for any $z> \hat{x}$.

Now we argue that on $S_{\tau_{\beta}}=\hat{x}$ we have $\lim\inf
W_{\tau_{\beta_n}} \geq\hat{x}$ almost surely. Coupled with the result
from the
previous paragraph we can then conclude that on
$W_{\tau_{\beta}}=\hat{x}$ we
have ${\tau_{\beta_n}} \rightarrow H_{\hat{x}} = \tau_\beta$.

Given $\delta$ and $\varepsilon< \hat{x} - \beta(\hat{x}-) - \delta$,
there exists
$N$ such that
for all $n>N$,
$\beta_n(\hat{x}-\varepsilon)<\beta(\hat{x}-) + \varepsilon<\hat
{x}-\delta$.
Then
\begin{eqnarray*}
\bigl( \omega\dvtx  W_{\tau_{\beta_n}}(\omega) < \hat{x}-\delta, S_{\tau
_{\beta}}(
\omega) = \hat{x} \bigr) & \subseteq& \bigl( \omega\dvtx  \inf\{W_u\dvtx
H_{\hat{x}-\varepsilon} \leq u \leq H_{\hat{x}} \} \leq\hat{x}-\delta
\bigr)
\\
&&{} \cup( \omega\dvtx  S_{\tau_{\beta_n}}< \hat{x} - \varepsilon, S_{\tau_{\beta
}} =
\hat{x}).
\end{eqnarray*}
By similar arguments to those in case 1 we can prove that the final
event has small probability.
Moreover, using that the fact that the probability that an event occurs
is
smaller than the expected number of times that it occurs,
\[
\Prob\bigl( \omega\dvtx  \inf\{W_u\dvtx  H_{\hat{x}-\varepsilon} \leq u \leq
H_{\hat{x}} \} \leq\hat{x}-\delta\bigr) \leq\int_{\hat{x}-\varepsilon
}^{\hat{x}}
\frac{dy}{y - (\hat
{x}-\delta)} = \ln\bigl(\delta/(\delta- \varepsilon)\bigr).
\]
By choosing $\varepsilon$ compared to $\delta$ this probability can be
made arbitrarily small.
\end{pf}

Note that if $\tau_{\beta_n} \rightarrow
\tau_\beta$ almost surely, then by the continuity of Brownian motion
$W_{\tau_{\beta_n}} \rightarrow W_{\tau_\beta}$ almost surely and
$\mu_n \Rightarrow\mu$.

We can summarize the results as follows:
%
%pr4.4 #&#
\begin{proposition}
\label{propmuUbeta}
Suppose that $(\mu_n)_{n \geq1}$ and $\mu$ are centered and that
$\E[|X_{\mu_n}|] \rightarrow\E[|X_\mu|]$.
Then the following are equivalent:
\begin{longlist}[(iii)]
\item[(i)] $\mu_n \Rightarrow\mu$ and $\E[|X_{\mu_n}|]
\rightarrow
\E[|X_{\mu}|];$
\item[(ii)]
$U_{n}(x) \rightarrow U_\mu(x)$ for each $x \in\R$;
\item[(iii)]
$\beta_n \rightarrow\beta$ at continuity points $s$ of $\beta$,
provided $s$ is less than or equal to the upper limit on the support of
$\mu$;
\item[(iv)]
$\tau_{\beta_n} \stackrel{\mathit{a.s.}}{\longrightarrow} \tau_\beta$;
\item[(v)] $W_{\tau_{\beta_n}} \stackrel{\mathit
{a.s.}}{\longrightarrow}
W_{\tau_\beta}$.
\end{longlist}
\end{proposition}
Now we want to prove a similar result for the Perkins embedding.
%
%le4.5 #&#
\begin{lemma}
\label{lemaconvergence}
Let $(\mu_n)_{n \geq1}$ be a sequence of centered probability
measures such that
$\mu_n \Rightarrow\mu$ and $\E[|X_n|] \rightarrow\E[|X_\mu|]$. Then
$a_n^\pm(x) \rightarrow a^{\pm}(x)$ at continuity points $x \in
(\check{x},\hat{x}) \setminus\{0\}$ of $a$. Moreover $\alpha^{\pm}_n(x)
\rightarrow\alpha^{\pm}(x)$ at nonzero continuity points $\check
{x}<x<\hat{x}$ of
$\mu$.
\end{lemma}
\begin{pf}
We prove the result for $(a_n^+,a^+)$, the other case being similar. The
extension from $a^{\pm}$ to $\alpha^{\pm}$ follows as in
Corollary~\ref{corbtobeta}.

Again we have that $x<0$ is a discontinuity point of $a^+$ if and only
if
there
is an atom of
$\mu$ at $x$. Suppose that $x$ is not an atom of $\mu$. Then, recall
(\ref{eqalphaplussolves}), $a^+(x)$ is the
unique solution in $z$ of $P(x)+P'(x)(z-x)=C(z)$. Moreover, for any
$\tilde{a}_n(x) \in(a_n^+(x+),a_n^+(x-))$,
\begin{eqnarray*}
P_n(x)+P'_n(x+) \bigl(
\tilde{a}_n(x)-x\bigr) &\geq& C_n\bigl(
\tilde{a}_n(x)\bigr),
\\
P_n(x)+P'_n(x-) \bigl(
\tilde{a}_n(x)-x\bigr) &\leq& C_n\bigl(
\tilde{a}_n(x)\bigr).
\end{eqnarray*}
Suppose $a^+_n(x) \rightarrow\gamma$ (down a subsequence if
necessary). Then since
$P'_n(x \pm) \rightarrow P'(x)$,
%and $P_n'(x-) \rightarrow P(x)$,
%
\[
P(x)+P'(x) (\gamma-x) \geq C(\gamma) \geq P(x)+P'(x)
(\gamma-x).
\]
Hence $\gamma=a^+(x)$ and $a_n^+(x) \rightarrow a(x)$.
\end{pf}

%alternative
%Suppose $C_n(x) \rightarrow C(x)$ pointwise.
%Then
%$a_n^\pm(x) \rightarrow a^{\pm}(x)$ at continuity points $x \in
%(\check{x},\hat{x}) \setminus\{0\}$ of $a$.
%
%Recall the definitions of $\theta= \theta_{C^{-1}}$ and $\theta_n =
%
%%If $y_n \in\partial C_n^{-1} (C_n(x))$ and if $ \partial C^{-1}$ is
%%single valued at $C(x)$ and takes the value $y$, say, then
%%\[(C^{-1}_n)^c(y_n) \rightarrow(C^{-1})^c(y)\]
%For $x>0$ %Figure~\ref{Fig2} demonstrates that
%$a^-_n(x)$ is the
%solution to the
%equation
%at least when $(C^{-1}_n)'(C_n(x))$ exists, see the figure below.
%When $\partial C_n^{-1}(C_n(x))$ is multivalued, we can extend the
%convention used for $b_n$ to define $a^-_n(x)$ to be the unique
%solution $a=a^-_n(x)$ of
%
%For $y<0$ define $L^{y}_n$ via $L^y_n(a) = a - P_n(a) y$ and
%$L^y$ via $L^y(a) =
%a - P(a)y$. Then, for example, $L^{\theta_n(x)}_n(a^{-}_n(x)) =
%b_n(x)$. Since $P_n$
%is continuous and decreasing in $a$, $L^y_n(a)$ is continuous and
%strictly
%increasing in $a$
%with continuous inverse
%$\ell^y_n$ say. Similarly, let $\ell^y$ denote the inverse to $L^y$.
%We have that if $y_n \rightarrow y$, then
%$L^{y_n}_n(a) \rightarrow L^y(a)$ pointwise, and
%$\ell^{y_n}_n(s) \rightarrow\ell^y(s)$. Finally, if $(C^{-1})'(C(x))$
%exists (and equals $\theta(x)$) then
%) = a^-(x). \]
%
%The proof for $a^+(x)$ is similar, but uses the reverse barycenter
%function, $\check{b}(x) = \E[X|X \leq x]$.
%
%pr4.6 #&#
\begin{proposition} \label{prop46}
\label{propmuUalpha}
Suppose that $(\mu_n)_{n \geq1}$ and $\mu$ are centered and that
$\E[|X_{\mu_n}|] \rightarrow\E[|X_\mu|]$.
\begin{enumerate}[(b)]
\item[(a)]
Suppose there exists an open interval $I$ containing $0$ such that
$\mu_n(I) = \mu(I)=0$.
Then the following are equivalent:
\begin{enumerate}[(iii)]
\item[(i)]
$\mu_n \Rightarrow\mu$ and $\E[|X_{\mu_n}|] \rightarrow\E
[|X_{\mu}|];$
\item[(ii)]
$U_{n}(x) \rightarrow U_\mu(x)$ for each $x \in\R$;
\item[(iii)]
$\alpha^{\pm}_n \rightarrow
\alpha^{\pm}$ at continuity points of $\alpha^{\pm}$ which lie
within the
range of the support of $\mu$;
\item[(iv)]
$\tau^{\rmP }_{\mu_n} \stackrel{\mathit{a.s.}}{\longrightarrow}
\tau^{\rmP }_{\mu}$;
\item[(v)] $W_{\tau^{\rmP }_{\mu_n}} \stackrel{\mathit
{a.s.}}{\longrightarrow}
W_{\tau^{\rmP }_{\mu}}$.
\end{enumerate}

\item[(b)]
More generally,
suppose $\mu_n \Rightarrow\mu$ and
$\E[|X_{\mu_n}|] \rightarrow\E[|X_{\mu}|]$.
Then,
$\alpha^{\pm}_n \rightarrow
\alpha^{\pm}$ at continuity points of $\alpha^{\pm}$ which lie
within the
range of the support of $\mu$.

Suppose further that
$\mu_{n}(\{0\}) \rightarrow\mu(\{0\})$. Then there exists a sequence of
Perkins embeddings of $\mu_n$ such that
$\tau_{\mu_n}^{\rmP }$ converges in probability to a Perkins embedding
$\tau^{\rmP }_\mu$
of
$\mu$.
In particular, if $Z_n$ converges in probability
to $Z$, then the Perkins embeddings
$(\tau^{\rmP ,Z_n}_{\mu_n})_{n \geq1}$ converge in probability to
the Perkins embedding $\tau_{\mu}^{\rmP ,Z}$ of
$\mu$.

Thus, if $\mu_n \Rightarrow\mu$,
$\E[|X_{\mu_n}|] \rightarrow\E[|X_{\mu}|]$ and $\mu_{n}(\{0\})
\rightarrow
\mu(\{0\})$, then if $(\tau^{\rmP ,Z_n}_{\mu_n})_{n \geq1}$ is a
sequence of
Perkins
embeddings of $(\mu_n)_{n \geq1}$, then there exists a subsequence
$n_k$ along
which
$\lim\tau^{\rmP ,Z_{n_k}}_{\mu_{n_k}}$ exists almost surely,
and is a Perkins embedding of $\mu$.
\end{enumerate}
\end{proposition}
\begin{pf}
For part (a) the equivalence of (i) and (ii) follows as
before. Lem\-ma~\ref{lemaconvergence} gives that (ii) implies
(iii). It follows from the pathwise construction of
$\tau_{\alpha_n}$ (and the existence of the interval $I$
which is not charged by $\mu_n$ so that $\tau^{\rmP }_{\mu_n} \equiv
\tau_{\alpha_n}$) that $\tau^{\rmP }_{\mu_n} \rightarrow\tau^{\rmP }_{\mu}$ almost
surely and hence
we have (iii) implies (iv). The continuity of Brownian
motion
allows us to deduce (v), from which (i) follows
immediately.
%; note that $\E[|X_{\mu_n}| \rightarrow\E[|X_\mu|]$ is one of the
%hypotheses of the Proposition.

For part (b) the statement about the convergence of $\alpha^{\pm}_n$ follows
as before. For the other results,
suppose
first that
$\mu(\{0\})=0$ and $\mu_n(\{0\})=0$ for
all
sufficiently large $n$.
Recall that $\tau_\alpha
= \inf\{ u\dvtx  W_u
\leq\alpha^+(S_u)$
or $W_u \geq\alpha^-(I_u) \}$ and for $\eta>0$ define the
stopping time
\[
\rho_{\alpha,\eta} = \tau_{\alpha_{\eta}}, %= \inf\{ u \geq H_{\eta}
\]
where $\alpha_{\eta}^+(s) = \min\{ \alpha^+(s), -\eta\}$,
$\alpha_{\eta}^-(i) = \max\{ \alpha^-(i), \eta\}$. Note that $\rho
_{\alpha,\eta}$ is the Perkins
embedding of a law which places no mass on $(-\eta,\eta)$.

We have that $\alpha_n \rightarrow\alpha$ at nonzero continuity points.
Let $\alpha_{n,\eta}^{\pm} =\break \mp\max\{ \mp\alpha_n^{\pm}(s),
\eta
\}$ and let $\rho_{\alpha_n,\eta}$ be the Perkins embedding for
$B_{\tau_{\alpha_n,\eta}}$.
Then $\alpha_{n,\eta}^{\pm} \rightarrow\alpha_{\eta}^{\pm}$ at
continuity points and by the
pathwise construction of $\rho_{\alpha_n,\eta}$, we have ${\rho_{\alpha
_n,\eta}} \rightarrow
{\rho_{\alpha,\eta}}$ almost surely.\vadjust{\goodbreak} In particular, given
$\delta, \varepsilon>0$ there exists $N_0$ such that for all $n \geq N_0$
\[
\Prob\bigl( | {\rho_{\alpha_n,\eta}} - {\rho_{\alpha,\eta}}| > \varepsilon\bigr) <
\delta/2.
\]

Note that on $|W_{\tau_{\alpha}}| > \eta$ we have
${\rho_{\alpha,\eta}} = {\tau_{\alpha}}$ with a similar statement
for $\alpha_n$. We can choose $\eta>0$ so that $\mu([-2 \eta,2 \eta
]) <
\delta/6$
and then $N_1$ so that for $n \geq N_1$, $\mu_n([- \eta, \eta]) <
\delta/3$.
Then
\begin{eqnarray*}
\bigl( | {\tau_{\alpha_n}} - {\tau_{\alpha}}| > \varepsilon\bigr) &\subseteq&\bigl(|
W_{\tau_{\alpha}} | \leq\eta\bigr) \cup\bigl(| W_{\tau_{\alpha_n}} |
\leq\eta\bigr)\\
&&{}\cup\bigl(|{
\tau_{\alpha_n}} - {\tau_{\alpha}} | > \varepsilon, | W_{\tau_{\alpha}} |
> \eta, | W_{\tau_{\alpha_n}} | > \eta\bigr)
\end{eqnarray*}
and the set $( | {\tau_{\alpha_n}} - {\tau_{\alpha}}| > \varepsilon)$
has probability at most $\delta$.

It follows that ${\tau_{\alpha_n}} \rightarrow{\tau_{\alpha}}$ in
probability, and hence that there is almost sure convergence down a
subsequence. Furthermore, down the same subsequence $W_{\tau_{\alpha_n}}
\rightarrow W_{\tau_{\alpha}}$ almost surely.

Now suppose that $\mu(\{0\})=0$ and that $\lim\mu_n(\{0\})=0$.
Recall the definition of $\mu_n^*$ as the measure $\mu_n$ with probability
mass at zero removed, and then rescaled to be a probability measure, and
note that $\alpha_{\mu_n^*} \equiv\alpha_{\mu_n}$. Then also $\mu_n^*
\Rightarrow
\mu$ and $U_{\mu_n^*} \rightarrow U_\mu$ pointwise.

Then, $\tau^{\rmP ,Z_n}_{\mu_n} = 0$ for $Z_n \leq\mu_{n}(\{0\})$ and
$\tau^{\rmP ,Z_n}_{\mu_n} = \tau_{\alpha_n}$ otherwise, so that $\tau
^{\rmP ,Z_n}_{\mu_n}
\rightarrow
\tau_\alpha$ in probability. Moreover, down a subsequence,
$\tau^{\rmP ,Z_n}_{\mu_n}
\rightarrow\tau_\alpha$ almost surely.

It remains to consider the case where $\mu(\{0\})>0$. For $\varepsilon<
1$, writing
$A_n = (Z_n \leq\mu_n(\{0\}), Z > \mu(\{0\}))$ and $B_n = (Z_n > \mu
_n(\{0\}), Z
\leq\mu(\{0\}))$,
\[
\bigl(\bigl| \tau^{\rmP ,Z_n}_{\mu_n} -\tau^{\rmP ,Z}_{\mu}
\bigr| > \varepsilon\bigr) \subseteq A_n \cup B_n \cup\bigl(
Z_n > \mu_n\bigl(\{0\}\bigr), Z > \mu\bigl(\{0\}\bigr),
|\tau_{\alpha_n} - \tau_{\alpha}|> \varepsilon\bigr)
\]
and $\tau^{\rmP ,Z_n}_{\mu_n}
\rightarrow\tau^{\rmP ,Z}_{\mu}$ in probability.
As before, there is almost sure convergence down a subsequence.
\end{pf}
%
%re4.7 #&#
\begin{remark}
One easy and natural way to guarantee that $Z_n \rightarrow Z$ is to
take $Z_n
= Z$ with probability one, or in other words to use the same independent
randomization
variable for each embedding.
\end{remark}
%
%re4.8 #&#
\begin{remark}
Suppose that $\mu$ is less than or equal to $\nu$ in convex order (we
write $\mu\leq_{\mathrm{cx}} \nu$). Then $U_{\mu} \leq U_\nu$. However, it does
not follow that $\beta_\mu\geq\beta_\nu$, and so it does not follow
that $\tau^{\mathrm{AY}}_{\mu} \leq\tau^{\mathrm{AY}}_{\nu}$. Similarly, we do not
have
that $|\alpha_\mu^{\pm}| \leq|\alpha_\nu^{\pm}|$ nor
$\tau^{\rmP }_{\mu}
\leq\tau^{\rmP }_{\nu}$.

Nonetheless, given $\mu$ it is possible to choose $\mu_n$ increasing to
$\mu$ in convex order and such that the barycenters are decreasing, and
hence\vspace*{1pt} the stopping times $\tau^{\mathrm{AY}}_{\mu_n}$ are monotonically
increasing and converge to $\mu$. This idea is used extensively
in Az\'ema and Yor~\cite{AzemaYor79a}, see also Revuz and
Yor~\cite{RevuzYor99}, Section VI.5, and also below in the proof of
Theorem~\ref{thmmain2}.

Similar remarks apply for the Perkins embedding.
\end{remark}
%
%ex4.9 #&#
\begin{example}
In Proposition~\ref{propmuUbeta} it does not hold that
$\beta_n(s) \rightarrow\beta(s)$ for $s$ beyond the upper limit on the
support of $\mu$.\vadjust{\goodbreak}

Suppose $\mu= \frac{1}{2}( \delta_{1} + \delta_{-1})$ and $\mu_n =
(1-n^{-2})
\frac{1}{2}( \delta_{1} + \delta_{-1}) +
n^{-2}\frac{1}{2} ( \delta_{n} + \delta_{-n})$. Then $U_{\mu}(0) = 1$
and
$U_{n}(0) = 1+ n^{-1}-n^{-2} \rightarrow1$.

We have $b_n$ is piecewise constant, and
$b_n(x)= 0$ for $x < -n$,
$b_n(x)= n/(2n^2-1)$ for $-n \leq x < -1$,
$b_n(x)= 1 + n^{-1}-n^{-2}$ for $-1 \leq x < 1$ and
$b_n(x)= n$ for $1 \leq x < n$.
Then $\beta_n(s) \rightarrow\beta_{\infty}(s)$ where
$\beta_{\infty}(s) = -1$ for $s \leq1$ and $\beta_{\infty}(s)=1$ for
$s>1$.
In contrast, $\beta(s)= -1$ for $s <1$ and $\beta(s)= s$ for $s \geq
1$.
\end{example}

%converge. For each $\omega$, there are subsequences $n(k)$ for which
%the
%path is stopped at either $\pm1$ (for the $\omega$ in the path, they
%are stopped at $+1$) and there are other subsequences $\tilde{n}(k)$
%for
%which the path is stopped when $|W|$ equals $1/\tilde{n}$, and so down
%these subsequences the stopping times converge to zero.}
%
%ex4.10 #&#
\begin{example}
If $\alpha_n \rightarrow\alpha_\mu$, but $U_n(0) \not\rightarrow
U_\mu(0)$, then in general
$\mu_n \not\Rightarrow\mu$.

Suppose $\mu= p( \delta_{1} + \delta_{-1}) + (1-2p) \delta_0$ and
$\mu_n = q( \delta_{1} + \delta_{-1}) + (1-2q) \delta_0$. Then
$\alpha_n
\equiv\alpha_\mu$ but $\mu_n \not\Rightarrow\mu$ unless $p=q$.
\end{example}
%
%ex4.11 #&#
\begin{example}
Suppose $\alpha_n \rightarrow\alpha_\mu$ at continuity points of
$\alpha_\mu$
and $U_n(0) \rightarrow U_\mu(0)$, but $\mu_n(\{0\})$ does not tend to
$\mu(\{0\})$. Then it does not follow that
$\tau_{\alpha_n}$ converges in probability, although even then we may
still have $\mu_n \Rightarrow\mu$.

Let
$\mu= \frac{1}{4}( \delta_{1} + \delta_{-1}) +
\frac{1}{2} \delta_0$,
and for $n>1$ let $\mu_n$ consist of masses of size
\[
\biggl\{ \frac{n+1}{4n}; \frac{1}{2}; \frac{n-1}{4n} \biggr\}
\]
at $\{ -1, 1/n, 1 \}$, respectively.
Then $\alpha^{\pm}(x) = \mp1$, $\alpha_n^{+}(x)= -1$ and $\alpha_n^-(x)
= 1/n$
for
$-1/n \leq x < 0$
and $\alpha_n^{-}(x)= 1$ for $x < - 1/n$. Further,
$\tau_\alpha= H_{\pm1}$ and
\[
\tau_{\alpha_n} = \cases{ H_{1/n}, &\quad if $H_{1/n}<H_{-1/n}$;
\cr
H_{-1}, &\quad if $H_{1/n}>H_{-1/n}$ and
$H_{-1}<H_{1}$;
\cr
H_{1}, &\quad if
$H_{1/n}>H_{-1/n}$ and $H_{1}<H_{-1}$.}
\]
Then, if $E_n$ is the event that $\tau_{\alpha_n} = H_{1/n}$, then
$\Prob(E_n)= 1/2$ and for $n>m$,
\[
\Prob(E_n|E_m) = \Prob(E_m|E_n)
= \Prob_{1/n}(H_{1/m} < H_{-1/m}) =
\frac{n +
m}{2n}.
\]
Hence $\Prob(E_n \cap E^c_m) = (n-m)/4n$ which does not tend to zero
as $n
\rightarrow\infty$ for fixed~$m$.
Hence
\[
\Prob\bigl(| \tau_{\alpha_n} - \tau_{\alpha_m}| > \varepsilon\bigr) \geq\Prob
\bigl( H_{\pm1} - H_{\pm1/2} > \varepsilon, E_n \cap
E_m^c\bigr) \not\rightarrow0
\]
and the
sequence $(\tau_{\alpha_n})_{n \geq1}$ is not Cauchy in probability.

%Suppose $\mu= \frac{1}{4}( \delta_{1} + \delta_{-1}) +
%and $\mu_n = \frac{1}{4}( \delta_{1} + \delta_{-1}) + \frac{1}{4}(
%Then $\alpha^{\pm}(x) = \mp1$ and $\alpha^{\pm}_n(x)= \mp1/n$ for
%$0< |x|
%1/(3n)$ and $\alpha^{\pm}_n(x)= \mp1$ for $|x|> 1/(3n)$. Further,
%$\tau_\alpha= \inf\{ u\dvtx  |W_u| = 1 \}$ and
% H_{-1/n} & \mbox{if $H_{-1/n}<H_{1/3n}$}; \\
% H_{1/n} & \mbox{if $H_{1/n}<H_{-1/3n}$}; \\
% H_{-1} & \mbox{if $H_{-1/n}>H_{1/3n}$, $H_{1/n}>H_{-1/3n}$ and
%$H_{-1}<H_{1}$}; \\
% H_{1} & \mbox{if $H_{-1/n}>H_{1/3n}$, $H_{1/n}>H_{-1/3n}$ and
%$H_{1}<H_{-1}$}. \end{array} \]
%Then, for almost every $\omega$, $\tau_{\alpha_n}(\omega)$ fails to
%converge, and there is both a subsequence converging to $0$, and
%another
%subsequence converging to $({H_{-1} \wedge H_1})(\omega)$.
\end{example}
%
%ex4.12 #&#
\begin{example}
Suppose $\alpha_n \rightarrow\alpha_\mu$ at continuity points of
$\alpha_\mu$ and $U_n(0) \rightarrow U_\mu(0)$ and
$\mu_{n}(\{0\})=0=\mu(\{0\})$. If there is no interval $I$ containing
$0$ on
which $\mu_n(I)=0=\mu(I)$, then it need not follow that $\tau_{\alpha_n}
\rightarrow\tau_\alpha$ almost surely, although there is convergence in
probability by Proposition~\ref{prop46}(b).

Let $\mu= U \{ -1, +1 \}$ and for $n>2$ let $\mu_n$ consist of masses
of size
\[
\biggl\{ \frac{n(1+2^{-n})}{2(1+n)}; \frac{1}{1+n}; \frac
{n(1-2^{-n})}{2(1+n)} \biggr\}
\]
at $\{ -1, n2^{-n}, 1 \}$, respectively.
Then $\alpha^{\pm}(x) = \mp1$, $\alpha_n^{+}(x)= -1$ and $\alpha_n^-(x)
= n2^{-n}$
for
$-2^{-n} \leq x < 0$
and $\alpha_n^{-}(x)= 1$ for $x < - 2^{-n}$. Further,
$\tau_\alpha= H_{\pm1}$ and
\[
\tau_{\alpha_n} = \cases{ H_{n2^{-n}}, &\quad if $H_{n2^{-n}}<H_{-2^{-n}}$;
\cr
H_{-1}, &\quad if $H_{n2^{-n}}>H_{-2^{-n}}$ and
$H_{-1}<H_{1}$;
\cr
H_{1}, &\quad if
$H_{n2^{-n}}>H_{-2^{-n}}$ and $H_{1}<H_{-1}$.}
\]
Then, if $E_n$ is the event that $\tau_{\alpha_n} \neq\tau_\alpha
$, then
$\Prob(E_n)= 1/(n+1)$ and for $n>m$,
\begin{eqnarray*}
\Prob(E_m \cap E_n) &=& \Prob(E_n)
\Prob(E_m |E_n) = \frac{1}{(1+n)} \frac{n2^{-n} + 2^{-m}}{m2^{-m}+2^{-m}}
%= \frac{1}{(1+n)} \frac{1}{(1+m)} (1+n2^{m-n})
\\
&=& \Prob(E_n) \Prob(E_m) +
\frac{n2^{m-n}}{(1+n)(1+m)}.
\end{eqnarray*}
Then by the Kochen--Stone lemma (Durrett~\cite{Durrett91}, Exercise 1.6.19),
$E_n$ happens infinitely often, almost surely. In particular,
almost surely, $\tau^{\rmP }_{\mu_n}$ does not converge.
\end{example}

%s5 #&#
\section{Objective functions as terminal values}
\label{secter}

Our goal is to prove that for a suitable class of bivariate functions
$F(w,s)$, the Az\'ema--Yor and
Perkins embeddings, which are well known to maximize and minimize $\E
[F(W_\tau,S_\tau)]$ in the special
case where $F$ does not depend on $w$ and $F$ is increasing in $s$,
continue to optimize this quantity
even when there is nontrivial dependence on $w$.
%
%as5.1 #&#
\begin{assumption}
\label{assF}
Throughout we assume that $F\dvtx  \{(w,s) \in\R\times\R_+;\break w \leq s \}
\mapsto\R_+$ is a continuous function and hence is bounded on compact
sets. We further assume that the partial derivative $F_s$ exists and is
continuous.
\end{assumption}

We are interested in functions $F$ which are monotonic in the following
sense (note that in our terminology increasing is a synonym for nondecreasing).
%
%de5.2 #&#
\begin{definition} \label{defFMON}
$F$ satisfies F-MON$\uparrow$ or F-MON$\downarrow$ if:
\begin{enumerate}[F-MON$\downarrow$]
\item[F-MON$\uparrow$]
$F_s(w,s)/(s-w)$ is monotonic increasing in $w$.
\item[F-MON$\downarrow$]
$F_s(w,s)/(s-w)$ is monotonic decreasing in $w$.
\end{enumerate}
\end{definition}

For $r \leq\hat{x} \leq\infty$ and $\eta\in
\{\beta,\alpha^+\}$ define
\[
\lambda_\eta(r) = \frac{F_s(\eta(r),r)}{r - \eta(r)},
\]
$\Lambda_\eta(s)=\int_0^s \lambda_\eta(r) \,dr$ and
$\Lambda^{(1)}_\eta(s)=\int_0^s r \lambda_\eta(r) \,dr$. Set
$\bar{\Lambda}_\eta= \sup_{s < \hat{x}} |\Lambda_{\eta}(s)|$.
Define $\Phi_\eta(w,s) = \int_0^s \lambda_\eta(r)(r-w) \,dr$;\vadjust{\goodbreak}
whence $\Phi_\eta(w,s) = \Lambda^{(1)}_\eta(s) - w \Lambda_\eta(s)$.
Finally, define $\xi_\beta(w)$ by
\[
\xi_\beta(w) = F\bigl(w, b(w)\bigr) - \Phi_{\beta}\bigl(w,
b(w)\bigr)
\]
and $\xi_{\alpha^+}(w)$ by
\[
\xi_{\alpha^+}(w) = F\bigl(w, \bar{a}(w)\bigr) - \Phi_{\alpha^+}
\bigl(w, \bar{a}(w)\bigr),
\]
where $\bar{a}(w)=w$ for $w \geq0$ and $\bar{a}(w)= a^+(w)$ for $w<0$.
Note that $\xi_\beta(w)$ [resp., $\xi_\alpha(w)$] does not depend
on the convention chosen for $b(w)$ [resp., $a^+(w)$].

%Suppose that $\mu$ has bounded support on $[\check{x},\hat{x}]$.
%Suppose that
%$\bar{\Lambda}_\eta< \infty$ for $\eta\in\{ \beta, \alpha^+ \}$.

%s5.1 #&#
\subsection{Target laws with bounded support}
\label{secbdd}

In this section we suppose $\mu$ has bounded support so that
$\check{x}$ and $\hat{x}$ are finite. This
assumption will be relaxed in the next section.
%
%th5.3 #&#
\begin{theorem}
\label{thmmain1a}
Suppose that $\mu$ has bounded support and that %Assumption
F-MON$\uparrow$ holds. Then
%Suppose $\bar{\Lambda}_\beta< \infty$. %and {\em integrability
%condition}. Then
%
%e5.1 #&#
%e5.2 #&#
\begin{eqnarray}
\label{eqthmAu} \sup_{\tau\in\sS_{\mathrm{UI}}(W,\mu)} \E\bigl[F(W_\tau,S_\tau)
\bigr] & = & \E\bigl[F(W_{\tau^{\mathrm{AY}}_\mu},S_{\tau^{\mathrm{AY}}_\mu})\bigr],
\\
%= \int F(w,b(w)) \mu(dw)
%Suppose $\bar{\Lambda}_\alpha< \infty$.
%and {\em integrability condition}.
%Then
\label{eqthmAl}
\inf_{\tau\in\sS_{\mathrm{UI}}(W,\mu)} \E\bigl[F(W_\tau,S_\tau)\bigr] & = & \E
\bigl[F(W_{\tau^{\rmP }_\mu},S_{\tau^{\rmP }_\mu})\bigr]. %=\int F(w,\bar{a}(w))
\end{eqnarray}
\end{theorem}

%If $F$ is nondecreasing in its second argument then \( \inf_{\tau\in
%To see this note that if $\tau$ embeds $\mu$ but is not
%uniformly integrable,
%then there exists $\sigma\in\sS_{\mathrm{UI}}(W,\mu)$ with $\sigma\leq\tau$
%and then $F(W_\tau,S_\tau) \geq F(W_\sigma,S_\sigma)$ almost surely.
%
%re5.4 #&#
\begin{remark}
\label{rembdef}
In the case where $\mu$ has no atoms
(so that the $\armin$ in (\ref{eqnbetadef}) is strictly increasing and
$\E[X|X \geq x] = \E[X|X>x]$),
then we can write
%
%e5.3 #&#
\begin{equation}
\label{eqnFint} \E\bigl[F(W_{\tau_\beta},S_{\tau_\beta})\bigr] = \int
_{\R} F\bigl(w,b_{\mu}(w)\bigr) \mu(dw).
\end{equation}
This formula need not hold if $\mu$ has atoms.

In cases where $\mu$ has a strictly positive density $\rho$ on
$(\check{x},\hat{x})$
and $\beta$
is differentiable,
the expression in (\ref{eqnFint}) can be rewritten as
%
%e5.4 #&#
\begin{eqnarray}
\label{eqnFint2} \E\bigl[F(W_{\tau_\beta},S_{\tau_\beta})\bigr] &=& \int
_{\R} F\bigl(\beta(s),s\bigr) \Prob(S_{\tau_\beta} \in ds) \nonumber\\[-8pt]\\[-8pt]
&=&
\int_{\R} F\bigl(\beta(s),s\bigr) \rho\bigl(\beta(s)\bigr)
\beta'(s) \,ds,\nonumber
\end{eqnarray}
where we use the fact that in the atom-free case
\[
\mu\bigl(\bigl[\beta(s),\infty\bigr)\bigr) = \Prob\bigl(W_{\tau_\beta} \geq\beta(s)
\bigr) = \Prob(S_{\tau_\beta} \geq s).
\]
A similar remark applies to
$\E[F(W_{\tau^{\rmP }_\mu},S_{\tau^{\rmP }_\mu})] =\int_{\R} F(w,\bar{a}(w))
\mu(dw)$.
%and the atom free case we have
%+ \int_{\R_+} F(\alpha(s),s)
\end{remark}
%
%re5.5 #&#
\begin{remark}
The requirement that the infimum in (\ref{eqthmAl}) is taken over
$\tau\in
\sS_{\mathrm{UI}}(W,\mu)$ (and not over all embeddings) is necessary, as can
be seen by
considering $F(w,s) = -(s-w)^3$. However, if we restrict attention to
functions $F$
which are increasing in $s$, then we may replace the infimum in
(\ref{eqthmAl}) with an infimum over all embeddings.\vadjust{\goodbreak}
\end{remark}

The key to the proof of the theorem is the following lemma.
%
%le5.6 #&#
\begin{lemma}
\label{lemmonup}
Suppose $F$ satisfies F-MON$\uparrow$. Then, for all $w \leq s$
\[
\xi_{\alpha^+}(w) + \Phi_{\alpha^+}(w,s) \leq F(w,s) \leq
\xi_{\beta}(w) + \Phi_{\beta}(w,s)
\]
with equality on the left at $w=s$ and $w=\alpha^+(w)$ and equality on
the
right at $w = \beta(s)$.
\end{lemma}
\begin{pf}
For $\eta\in\{ \beta, \alpha^+ \}$ define
%
%e5.5 #&#
\begin{equation}
\label{eqinlagrange} L_{\eta} (w,s)= \biggl[F(w,s)-
\xi_\eta(w)- \int_0^s
\lambda_\eta(r) (r-w)\,dr \biggr].
\end{equation}

We will show that
$L_{\alpha^+}(w,s) \geq0$ with equality at $w=s$ and $w=\alpha^+(s)$,
and
$L_\beta(w,s) \leq0$ with equality at $w=\beta(s)$.

Consider the latter inequality first:
\begin{eqnarray*}
L_\beta(w,s) &=& F(w,s)-\xi_\beta(w)- \int
_0^s \lambda_\beta(r) (r-w)\,dz
\\
&=&F(w,s)-F\bigl(w,b(w)\bigr)+\int_0^{b(w)} \,dr
F_s\bigl(\beta(r),r\bigr) \frac{r-w}{r-\beta(r)} \\
&&{}- \int
_0^s dr F_s\bigl(\beta(r),r\bigr)
\frac{r-w}{r-\beta(r)}
\\
&=& \int_{b(w)}^s \biggl\{ \frac{F_s(w,r)}{r-w} -
\frac{F_s(\beta(r),r)}{r-\beta(r)} \biggr\} (r-w) \,dr.
\end{eqnarray*}
If $b(w) < r < s$, then since $\beta$ is increasing,
$w<\beta(r)$ and by F-MON$\uparrow$ the integrand is negative. If
$s<r<b(w)$, then
$w>b(r)$ and the integrand is positive. Thus $L_\beta(w,s) \leq0$
as required. Clearly, there is equality at $s=b(w)$.

For $L_{\alpha^+}$ a similar calculation to the one above shows that
\[
L_{\alpha^+}(w,s) = \int_{\bar{a}(w)}^s \biggl\{
\frac
{F_s(w,r)}{r-w} - \frac{F_s(\alpha^+(r),r)}{r-\alpha^+(r)} \biggr\}
(r-w) \,dr.
\]
To see that $L_{\alpha^+}(w,s) \geq0$, consider $w\geq0$ and $w<0$
separately.
For $w \geq0$, $\bar{a}(w)=w$ and for $w<r<s$, $\alpha^+(r) \leq
\alpha^+(w) \leq w$ so
that the integrand is positive and $L_{\alpha^+}(w,s) \geq0$.
For $w<0$, $\bar{a}(w)=a(w)$, and then if $a(w)<r<s$, we have
$w>\alpha^+(r)$ and the integrand is positive. Otherwise if
$s<r<a(w)$, $w<\alpha^+(r)$ and
the integrand is negative. In either case, allowing for the limits on
the integral, $L_{\alpha^+}(w,s) \geq0$.
Equality holds at $w=s$ and $w=\alpha^+(s)$.
\end{pf}
%
%re5.7 #&#
\begin{remark} \label{ridearelated}
Essentially, the idea behind Lemma~\ref{lemmonup} and the proof of Theorem
\ref{thmmain1a} is to
interpret the embedding property and
Doob's (in)-equality for the martingale $W$ as linear constraints on
the
possible joint laws of $(W_\tau,S_\tau)$,\vadjust{\goodbreak} with associated Lagrange
multipliers.
Thus, if the joint law is given by $\nu(dw,ds)$, then
$\int_{s \geq r} (w-r) \nu(dw,ds) = 0$ (which is equivalent to (3.2) in
Rogers~\cite{Rogers93}).
There is an identity of this form for each $r$ and
when they are integrated against a family of Lagrange multipliers
$\lambda_\eta(r)$ we
obtain
\[
0 = \int_{0}^{\infty} \lambda_\eta(r)
\int_{s \geq r} (w-r) \nu(dw,ds) = \int\nu(dw,ds) \int
_{0 \leq r \leq s} \lambda_\eta(r) (w-r) \,dr.
\]
The integrand of this last expression appears as the last term in
(\ref{eqinlagrange}).
\end{remark}

It remains to prove Theorem~\ref{thmmain1a}.
The main idea for the proof of the theorem is that provided that
$\bar{\Lambda}_\beta$ and $\bar{\Lambda}_{\alpha^+}$ are finite,
then for $\tau\in
\sS_{\mathrm{UI}}(W,\mu)$ both
$(\Phi_{\alpha^+}(W^\tau_t,S^\tau_t))_{t \geq0}$ and
$(\Phi_{\beta}(W^\tau_t,S^\tau_t))_{t \geq0}$ are
uniformly integrable martingales.
[By It\^{o}'s formula, $d \Phi_{\eta}(W_t,S_t) = - \Lambda_\eta(S_t)
\,dW_t$ since the finite variation term involves the product $(S_t - W_t)
\,dS_t$ and when $S$ is increasing we must also have $S_t - W_t=0$.]
It follows that $\E[\Phi_{\beta}(W_\tau,S_\tau)]=0$ and
\[
\E\bigl[\xi_{\alpha^+}(W_\tau)\bigr] \leq\E\bigl[F(W_\tau,S_\tau)
\bigr] \leq\E\bigl[\xi_{\beta}(W_\tau)\bigr],
\]
which, given the forms of $\xi_{\alpha}$ and $\xi_{\beta}$ leads to the
first result given in the \hyperref[sec-intro]{Introduction}.
%
%re5.8 #&#
\begin{remark} \label{remarkazemayor}
The processes $(\Phi_{\alpha^+}(W^\tau_t,S^\tau_t))_{t \geq0}$ and
$(\Phi_{\beta}(W^\tau_t,S^\tau_t))_{t \geq0}$ belong to the class of
Az\'{e}ma--Yor martingales. A martingale $M=(M_t)_{t \geq0}$ is an
Az\'{e}ma--Yor martingale if $M_t=G(S^X_t)-(S^X_t-X_t)g(S_t)$
for $X$ a martingale and $G'=g$; see~\cite{AzemaYor79a}.
\end{remark}
%
%re5.9 #&#
\begin{remark}
An alternative derivation of (the right inequality of) Lem\-ma~\ref
{lemmonup} is to look for pathwise
inequalities
$F(W_t,S_t) \leq\xi(W_t) + M_t$ such that $M_t$ is a Markovian
function of $W_t$ and $S_t$ and such that
there is equality at $S_t = b(W_t)$.

If $M_t = \Phi(W_t,S_t)$ and $\Phi$ is appropriately differentiable,
then $M$ must be an Az\'ema--Yor
martingale $\Phi(W_t,S_t) = -H(S_t) + H'(S_t)(S_t-W_t)$ for some~$H$.
Further, if there is to be equality
at
$s=b(w)$, then we must have $\xi(w) = F(w,b(w)) - \Phi(w,b(w))$. Then
we want conditions on $F$ such that
there is an inequality $F(w,s) \leq\xi(w) + \Phi(w,s)$, or equivalently
\begin{eqnarray*}
\int_{b(w)}^{s} F_s(w,r) \,dr &=& F(w,s) -
F\bigl(w,b(w)\bigr) \\
& \leq& \Phi(w,s) - \Phi\bigl(w,b(w)\bigr) = \int
_{b(w)}^{s} \Phi_s(w,r) \,dr
\\
& = & \int_{b(w)}^{s} H''(r)
(r-w) \,dr.
\end{eqnarray*}
From this it follows that a sufficient condition is $F_s(w,r) \leq
H''(r)(r-w)$ for $r > b(w)$ and the
reverse inequality for $r<b(w)$, which holds if F-MON$\uparrow$ holds
and $H''(s) =
F_s(\beta(s),s)/(s-\beta(s))$.\vadjust{\goodbreak}
\end{remark}
\begin{pf*}{Proof of Theorem~\ref{thmmain1a}}
Consider first the bound associated with the Az\'{e}ma--Yor embedding.
$\bar{\Lambda}_\beta$ depends on the
combination of $\mu$ and $F$.

Suppose that $\mu$ has an atom at $\hat{x}$.
By Lemma~\ref{lemintegrabilitybeta}
$(r-\beta(r))^{-1}$ is integrable near zero so that if $\mu$ has an atom
at $\hat{x}$, then $r-\beta(r)$ is bounded below for $r < \hat{x}$ and
$\bar{\Lambda}_\beta< \infty$.
Since $\tau\in\sS_{\mathrm{UI}}(W,\mu)$ implies $(W^\tau_t)_{t \geq0}$ is
bounded, and since %by Assumption~\ref{assbdd}
$\Lambda_\beta(s)$ and
$\Lambda^{(1)}(s)$ are bounded,
we have that $\Phi_{\beta}(W_{t}^\tau, S_{t}^\tau)$ is a
bounded local martingale and hence
$\E[ \Phi_{\beta}(W_{t}^\tau, S_{t}^\tau)] = 0$, which can be
re-expressed
as
$\E[ \Lambda_\beta^{(1)}(S_\tau)] = \E[W_\tau
\Lambda_\beta(S_\tau)]$.
In view of Lemma~\ref{lemmonup} we have
%
%e5.6 #&#
\begin{equation}
\label{pathwiseineq} F(W_{\tau}, S_{\tau}) \leq
\xi_\beta(W_{\tau}) + \Phi_{\beta}(W_{\tau},
S_{\tau}).
\end{equation}
Thus
\[
\E\bigl[F(W_\tau,S_\tau)\bigr] \leq\int\xi_\beta(w)
\mu(dw).
\]
Note that for $\tau= \tau_{\beta}$, we have equality in
(\ref{pathwiseineq}) and hence equality in this last expression.

Now suppose there is no atom at $\hat{x}$. Fix $\tau\in\sS_{\mathrm{UI}}(W,\mu)$
and let
$\sigma_n = \tau\wedge H_{\check{x}-1/n}$ and $\mu_n = \sL
(W_{\sigma_n})$. Then
$U_{\mu_n} \rightarrow U_\mu$ for each $x$ and by bounded convergence
we have both
\[
\E\bigl[F(W_\tau,S_\tau)\bigr] = \E\bigl[\lim
F(W_{\sigma_n},S_{\sigma_n})\bigr] = \lim\E\bigl[F(W_{\sigma
_n},S_{\sigma_n})
\bigr]
\]
and
\[
\E\bigl[F(W_{\tau^{\mathrm{AY}}_{\mu}},S_{\tau^{\mathrm{AY}}_{\mu}})\bigr] = \E\bigl[\lim
F(W_{\tau^{\mathrm{AY}}_{\mu_n}},S_{\tau^{\mathrm{AY}}_{\mu_n}})\bigr] = \lim\E\bigl
[F(W_{\tau^{\mathrm{AY}}_{\mu_n}},S_{\tau^{\mathrm{AY}}_{\mu_n}})
\bigr].
\]
The result follows from the previous case on comparing $\sigma_n$ with
$\tau^{\mathrm{AY}}_{\mu_n}$.

The proof of (\ref{eqthmAl}) is identical except that there is no
need to
separate the case where there is an atom at $\hat{x}$, since by
Lemma~\ref{lemintegrabilityalpha} $(r-\alpha^+(r))^{-1}$ is
integrable near zero and hence the fact
that $\mu$ has bounded support is sufficient for
$\bar{\Lambda}_{\alpha^+}<\infty$.
\end{pf*}

There are a parallel pair of results based on F-MON$\downarrow$, the
proofs of which are very similar.
%
%le5.10 #&#
\begin{lemma}
\label{lemmondown}
Suppose $F$ satisfies F-MON$\downarrow$. Then, for all $w \leq s$
\[
\xi_{\beta}(w) + \Phi_{\beta}(w,s) \leq F(w,s) \leq
\xi_{\alpha^+}(w) + \Phi_{\alpha^+}(w,s)
\]
with equality on the right at $s=w$ and $s=a(w)$ and equality on the
left at $s = b(w)$.
\end{lemma}
%
%th5.11 #&#
\begin{theorem}
\label{thmmain1b}
Suppose F-MON$\downarrow$ holds. %and Assumption~\ref{assbdd} hold.
Then
%Suppose $\bar{\Lambda}_\beta< \infty$. % and
%{\em integrability condition}. Then
%
\begin{eqnarray*}
\inf_{\tau\in\sS(W,\mu)} \E\bigl[F(W_\tau,S_\tau)\bigr] &=& \E
\bigl[F(W_{\tau_\mu^{\mathrm{AY}}},S_{\tau_\mu^{\mathrm{AY}}})\bigr],
\\
%
%= \int\xi_{\beta}(w) \mu(dw) < \infty\]
%Suppose $\bar{\Lambda}_\alpha< \infty$.
%and {\em integrability condition}.
%Then
%
\sup_{\tau\in\sS_{\mathrm{UI}}(W,\mu)} \E\bigl[F(W_\tau,S_\tau)\bigr] &=& \E
\bigl[F(W_{\tau_\mu^{\rmP }},S_{\tau_\mu^{\rmP }})\bigr].
\end{eqnarray*}
%
%=\int\xi_{\alpha}(w) \mu(dw) < \infty\]
\end{theorem}
%
%ex5.12 #&#
\begin{example} \label{expowerunif}
Suppose $\mu= U[-1,1]$ and
$F(w,s)=(s-w)^c$ for $c > -1$ (with $c \neq0$). Then for $c\geq2$
F-MON$\downarrow$
holds, for $0<c\leq2$ F-MON$\uparrow$ holds and for $-1<c<0$,
F-MON$\downarrow$
holds again.

Write $B^{\mathrm{AY}}$ and $B^{\rmP }$ for the bounds associated with the
Az\'ema--Yor and Perkins embeddings.

Recall the expressions for $\beta$ and $\alpha$ from
Examples~\ref{excalcuniform} and~\ref{excalcuniform2}.

For the Az\'{e}ma--Yor embedding, $\beta(s)=2s-1$ and the law of the
$S_{\tau_\beta}$ is a
uniform on $[0,1]$. The associated bound (as a function of the parameter
$c$) is
given by
\begin{eqnarray*}
B^{\mathrm{AY}}(c) &=& \E\bigl[F(W_{\tau^{\mathrm{AY}}_\mu},S_{\tau^{\mathrm{AY}}_mu})\bigr] =\int
_{-1}^1 \bigl(b(w)-w\bigr)^c
\,\frac{dw}{2} =\int_0^1 \bigl(s-\beta(s)
\bigr)^c \,ds \\
&=&\int_0^1
(1-s)^c \,ds = \frac{1}{c+1}.
\end{eqnarray*}

For the Perkins bound, note that for $c<0$, $F(s,s)= \infty$, and so
$B^{\rmP }(c)=0$. For $c>0$, $F(s,s)=0$ and using the substitution $w =
\alpha^+(s)=s - 2 \sqrt{s}$,
\begin{eqnarray*}
B^{\rmP }(c) &=& \E\bigl[F(W_{\tau^{\rmP }_\mu},S_{\tau^{\rmP }_{\mu}})\bigr] =
\int_{-1}^0 \bigl(a^+(w)-w\bigr)^c
\,\frac{dw}{2}
\\
%&=& \int_0^1 (s-\alpha^+(s))^c \Prob(S_{\tau_\mu^{\rmP }} \geq s,
%W_{\tau_\mu^{\rmP }} <s) \,ds \\
%&=& \int_0^1 (s-\alpha^+(s))^c \frac{\Prob(S_{\tau_\mu^{\rmP }} \geq s)
%ds}{s-\alpha^+(s)} \\ &=& \int_0^1 (2 \sqrt{s})^{c-1} (1-\sqrt{s}) \,ds
%
%
&=&
\frac{2^c}{(c+1)(c+2)}.
\end{eqnarray*}

Results for a range of $c$ are plotted in Figure~\ref{fig3}. Observe that for
$c=2$, $B^{\mathrm{AY}}(2)= B^{\rmP }(2)=1/3$ and all uniformly
integrable embeddings for the terminal law are consistent with the same
expected payoff. The reason for this will become clear in Section
\ref{sec-g} and will correspond to the choice $g \equiv1$.

%f3 #&#
\begin{figure}[b]

\includegraphics{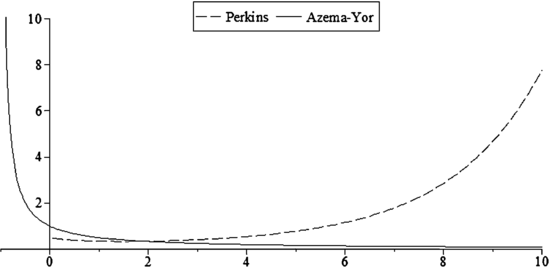}

\caption{All uniformly integrable embeddings have the same expected
value when $c=2$. Note the reversal of the bounds at $c=2$:
for $0<c<2$ Theorem \protect\ref{thmmain1a} applies while for $c>2$
Theorem \protect\ref{thmmain1b} applies.
For $c<0$, the Perkins bound is infinite and the Az\'{e}ma--Yor bound
is finite.
The Perkins bound as a function of $c$ is discontinuous at $c=0$.}\label{fig3}
\end{figure}

In fact Assumption~\ref{assF} is not satisfied for $-1<c<1$.
Nonetheless, for $c$ in
this range and
$\varepsilon>0$ we can let\vadjust{\goodbreak} $F_\varepsilon(w,s)= h_\varepsilon(s-w)$ where
$h_\varepsilon(x) = x^c$ for $x\geq\varepsilon$
and $h_\varepsilon(x) = \varepsilon^c + c \varepsilon^{c-1}(x-c)$ for
$x<\varepsilon$. Then $F_\varepsilon$ does satisfy
Assumption~\ref{assF}, and $F$ and $F_\varepsilon$ satisfy
F-MON$\uparrow$ or
F-MON$\downarrow$ together. Then arguments of Theorem \ref
{thmmain1a} provide the upper and lower bounds
for $F_\varepsilon$, and letting $\varepsilon\downarrow0$ we obtain the
pictured bounds for $F$.
\end{example}
%
%ex5.13 #&#
\begin{example} \label{expowerunif2}
Suppose again that $\mu=U[-1,1]$. Let $F(w,s)= \frac{(s-w)^2}{s^c}$.
Note that for each $c$ either F-MON$\uparrow$ or
F-MON$\downarrow$ (or both) holds, so that the Az\'ema--Yor and Perkins
embeddings give extremal values for $\E[F(W_\tau,S_\tau)]$.
%
%f4 #&#
\begin{figure}

\includegraphics{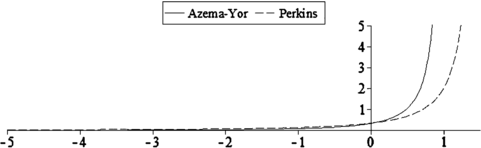}

\caption{For $1<c<3/2$ the Az\'ema--Yor upper bound is infinite while the
Perkins lower
bound is finite.}\label{fig4}\vspace*{-3pt}
\end{figure}
Consider the Az\'ema--Yor bound as a function of the
parameter $c$ (defined for $c<1$),
\begin{eqnarray*}
B^{\mathrm{AY}}(c)&=& \int_{-1}^1
\frac{(b(w)-w)^2}{b(w)^c} \,\frac{dw}{2} =\int_0^1
\frac{(s-\beta(s))^2}{s^c} \,ds =\int_0^1
\frac{(s-1)^2}{s^c} \,ds \\
&=&\frac{2}{(1-c)(2-c)(3-c)}.
\end{eqnarray*}
For the Perkins bound we have (for $c<3/2$)
\begin{eqnarray*}
B^{\rmP }(c) &=& \int_{-1}^0
\frac{(a^+(w) -w)^2}{a^+(w)} \,\frac{dw}{2}
\\
%&=& \int_0^1 \frac{(s-\alpha^+(s))^2}{s^c}\frac{\Prob(S_{\tau_\mu^{\rmP }}
%>s) \,ds}{s-\alpha^+(s)} \\
&=& \int_0^1
\frac{2 \sqrt{s}}{s^c} (1-\sqrt{s}) \,ds
\\
&=& \frac{1}{(3/2-c)(2-c)}.
\end{eqnarray*}

Observe that the expressions for $B^{\mathrm{AY}}(\cdot)$ and $B^{\mathrm{P}}(\cdot)$
co-incide at $c=0$ where both F-MON$\uparrow$ and F-MON$\downarrow$
hold. See Figure~\ref{fig4}.\vspace*{-2pt}
\end{example}

%s6 #&#
\section{General centered target measures}
\label{sec-gen}

%The goal of this section is to extend the results of
%Section~\ref{secbdd}
%to a wider class of target distributions, and simultaneously to
%relax the
%assumption that $\bar{\Lambda}_\eta$ is finite.

%Given $\tau\in\sS_{\mathrm{UI}}(W,\mu)$, let $\sigma_n = \tau\wedge H_{\pm
%n}$, $\mu_n = \sL(W_{\sigma_n})$ and define $\tau^{\mathrm{AY}}_{\mu_n}$ and
%$\tau^{\rmP }_{\mu_n}$ to be the
%Az\'ema--Yor and Perkins
%stopping times associated with $\mu_n$.
%
%th6.1 #&#
\begin{theorem}
\label{thmmain1c}
Fix $\tau\in\sS_{\mathrm{UI}}(W,\mu)$. Suppose, in addition to
Assumption~\ref{assF},
that $F \geq0$, that
%
%e6.1 #&#
\begin{equation}
\label{eqnuiF} \E\bigl[ F(W_{H_{\pm n}}, S_{H_{\pm n}}); \tau\geq
H_{\pm n} \bigr] \rightarrow0\vadjust{\goodbreak}
\end{equation}
and that if $(\mu_n)_{n \geq1}$ is any sequence of measures which is
increasing
in convex order for which $\mu_n \Rightarrow\mu$, $U_{\mu_n}(0)
\rightarrow
U_{\mu}(0)$
and $\mu_{n}(\{0\})
\rightarrow\mu(\{0\})$, then both
%
%e6.2 #&#
%e6.3 #&#
\begin{equation}
\label{eqnuib}
\E\bigl[ F(W_{\tau^{\mathrm{AY}}_{\mu_n}},
S_{\tau^{\mathrm{AY}}_{\mu _n}}) \bigr] \rightarrow \E\bigl[
F(W_{\tau^{\mathrm{AY}}_{\mu}}, S_{\tau^{\mathrm{AY}}_{\mu}}) \bigr ]
\end{equation}
and
\begin{equation}
\label{eqnuia}
\E\bigl[ F(W_{\tau^{\rmP }_{\mu_n}}, S_{\tau^{\rmP
}_{\mu_n}}) \bigr] \rightarrow \E\bigl[ F(W_{\tau^{\rmP }_{\mu}},
S_{\tau^{\rmP }_{\mu}}) \bigr].
\end{equation}
Then if F-MON$\uparrow$ holds,
\[
\E\bigl[ F(W_{\tau^{\rmP }_{\mu}}, S_{\tau^{\rmP }_{\mu}}) \bigr] \leq\E\bigl[
F(W_{\tau}, S_{\tau}) \bigr] \leq\E\bigl[ F(W_{\tau^{\mathrm{AY}}_{\mu}},
S_{\tau^{\mathrm{AY}}_{\mu}}) \bigr],
\]
whereas, if F-MON$\downarrow$ holds, then
\[
\E\bigl[ F(W_{\tau^{\mathrm{AY}}_{\mu}}, S_{\tau^{\mathrm{AY}}_{\mu}}) \bigr] \leq\E
\bigl[
F(W_{\tau}, S_{\tau}) \bigr] \leq\E\bigl[ F(W_{\tau^{\rmP }_{\mu}},
S_{\tau^{\rmP }_{\mu}}) \bigr].
\]
\end{theorem}
\begin{pf}
Suppose F-MON$\uparrow$ holds (the proof for F-MON$\downarrow$ is
similar).
Given $\tau\in\sS_{\mathrm{UI}}(W,\mu)$, let $\sigma_n = \tau\wedge H_{\pm
n}$, $\mu_n = \sL(W_{\sigma_n})$ and define $\tau^{\mathrm{AY}}_{\mu_n}$ and
$\tau^{\rmP }_{\mu_n}$ to be the
Az\'ema--Yor and Perkins
stopping times associated with $\mu_n$.

We have, using monotone convergence,
(\ref{eqnuiF}), Theorem~\ref{thmmain1a} and finally (\ref{eqnuib}),
\begin{eqnarray*}
\E\bigl[ F(W_{\tau}, S_{\tau}) \bigr] & = & \E\bigl[ \lim
F(W_{\sigma_n}, S_{\sigma_n}); \sigma_n = \tau\leq
H_{\pm n} \bigr]
\\
& = & \lim\E\bigl[ F(W_{\sigma_n}, S_{\sigma_n}) I_{\{ \tau\leq H_{\pm
n} \} }
\bigr]
\\
& = & \lim\E\bigl[ F(W_{\sigma_n}, S_{\sigma_n})I_{ \{ \tau
< H_{\pm n} \} } +
F(W_{H_{\pm n}}, S_{H_{\pm n}})I_{\{ \tau\geq H_{\pm n}
\} } \bigr]
\\
& = & \lim\E\bigl[ F(W_{{\sigma_n}}, S_{{\sigma_n}}) \bigr]
\\
& \leq& \lim\E\bigl[ F(W_{\tau^{\mathrm{AY}}_{\mu_n}}, S_{\tau^{\mathrm{AY}}_{\mu_n}})
\bigr] = \E\bigl[
F(W_{\tau^{\mathrm{AY}}_\mu}, S_{\tau^{\mathrm{AY}}_\mu}) \bigr].
\end{eqnarray*}

Similarly
\[
\lim\E\bigl[ F(W_{{\sigma_n}}, S_{{\sigma_n}}) \bigr] \geq\lim\E\bigl[
F(W_{\tau^{\rmP }_{\mu_n}}, S_{\tau^{\rmP }_{\mu_n}}) \bigr] = \E\bigl[ F(W_{\tau
^{\rmP }_{\mu}},
S_{\tau^{\rmP }_{\mu}}) \bigr].
\]
\upqed\end{pf}
%
%co6.2 #&#
\begin{corollary}
\label{corl2}
Suppose that $F(w,s) \leq A(1 + |w|^k +
s^k)$ for $k \geq1$ and that $\mu$ has finite $k+\varepsilon$ moment,
for some positive $\varepsilon$. Then the
hypotheses
(\ref{eqnuiF}), (\ref{eqnuib}) and
(\ref{eqnuia}) are all satisfied, and the conclusions of
Theorem~\ref{thmmain1c} hold.
\end{corollary}
\begin{pf}
By Doob's submartingale inequality for $(|W_{t \wedge\tau}|^{k +
\varepsilon})_{t \geq0}$, for any $\tau\in\sS_{\mathrm{UI}}(W, \mu)$,
\[
m^{k+\varepsilon} \Prob( \tau> H_{\pm m} ) < \E\bigl[|W_\tau|^{k+\varepsilon}
\bigr] < \infty.
\]
Then
\[
\E\bigl[ F(W_{H_{\pm n}}, S_{H_{\pm n}}); \tau\geq H_{\pm n}
\bigr] \leq A\bigl(1 + 2 n^k\bigr) \Prob( \tau> H_{\pm n} )
\rightarrow0.
\]

For (\ref{eqnuib}) %(the proof of (\ref{eqnuia}) is similar)
we have that $\tau_{\beta_n} \rightarrow\tau_{\beta}$ almost
surely. Moreover, since $\mu_n \leq_{\mathrm{cx}} \mu$ there exists a
stopping time ($\rho_n$ say) with $\rho_n \geq\tau_{\beta_n}$ and
$\rho_n \in\sS_{\mathrm{UI}}(W,\mu)$. For such a~$\rho_n$,
$\E[|W_{\rho_n}|^{k+\varepsilon}] = \int_{\R}|x|^{k+\varepsilon} \mu
(dx) <
\infty$ by hypothesis, and then (defining $W^*_t= \sup_{s \leq t}|W_s|$)
by Doob's\vspace*{1pt} $L^p$ inequality $\E[(W^*_{\rho_n})^{k + \varepsilon}] \leq D
< \infty$
for some constant $D$, independent of $n$.

Set\vspace*{1pt} $F_n = F(W_{\tau^{\mathrm{AY}}_{\mu_n}}, S_{\tau^{\mathrm{AY}}_{\mu_n}})$
and $F = F(W_{\tau^{\mathrm{AY}}_{\mu}}, S_{\tau^{\mathrm{AY}}_{\mu}})$, then $F_n
\rightarrow F$ almost surely. The goal is to show that $\E[F_n]
\rightarrow\E[F]$ which will follow if\break
$\sup_n \E[(F_n)^p] < \infty$, for then
$(F_n)_{n \geq1}$ is
uniformly integrable.
We have that if $|w| \leq x$ and $s \leq x$, then with $p = 1 +
k/\varepsilon$,
\[
F(w,s)^p \leq A^p\bigl(1 + 2x^k
\bigr)^p \leq A^p 3^p \bigl(1+
x^{kp}\bigr).
\]
Hence
\[
\E\bigl[ F_n^p \bigr] \leq A^p
3^p \bigl( 1 + \E\bigl[\bigl(W^*_{\tau_n}
\bigr)^{kp}\bigr]\bigr) \leq A^p 3^p \bigl( 1
+ \E\bigl[\bigl(W^*_{\rho_n}\bigr)^{kp}\bigr]\bigr) \leq
A^p 3^p ( 1 + D) < \infty.
\]

For (\ref{eqnuia}), consider a subsequence $n(k)$. Then down a further
subsequence $\tau^{\rmP }_{\mu_n} \rightarrow\tau^{\rmP }_{\mu}$ almost
surely and down this subsequence (\ref{eqnuia}) holds by identical
arguments as in the case for the Az\'ema--Yor embedding. Hence
(\ref{eqnuia}) holds.
\end{pf}

%s7 #&#
\section{Objective functions as running costs}
\label{sec-g}

Our original aim in studying functions $F(w,s)$ was as an aid in
the
analysis of the expected values of integrals of the form $\int_{0}^\tau
g(S_t) \,dt$.
Motivated by the variance swap problem in mathematical finance we asked:
\begin{quote}
Given $g$ and $\mu$, what is the range of possible values of
$\E[\int_0^\tau g(S_u) \,du]$ over embeddings $\tau$ of $\mu$ in Brownian
motion?
\end{quote}
Our aim is to reduce this problem to the case previously considered, but
to use the extra structure to prove more powerful results under weaker
hypotheses.

The expected value of $\int_0^\tau g(S_u) \,du$ is intimately related
to the value of $\E[ G(W_\tau, S_\tau) ]$
where $G(w,s)=(s-w)^2 g(s)$. Indeed, if $g$ is continuously
differentiable,
then by
It\^{o}'s lemma,
%
%e7.1 #&#
\begin{equation}
\label{eqnitoG} G(W_{\tau},S_{\tau}) = G(0,0) + \int
_{0}^{\tau} g(S_u) \,du - \int
_0^{\tau} 2(S_u - W_u)
g(S_u) \,dW_u,
\end{equation}
so that if $g(0)$ is finite [and then $G(0,0)=0$], and if
\[
\biggl(\int_0^{\tau\wedge t} 2(S_u - W_u) g(S_u) \,dW_u\biggr)_{t \geq0}
\]
is a uniformly integrable martingale, then $\E[ \int_0^\tau
g(S_u) \,du ] = \E[ G(W_\tau, S_\tau) ]$.

If $g$ is increasing (resp., decreasing), then $G$ satisfies
G-MON$\downarrow$ (resp., G-MON$\uparrow$), and we can apply
the results of previous sections to deduce that the Az\'{e}ma--Yor
and Perkins solutions give bounds
$\E[ \int_0^\tau g(S_u) \,du ]$ over embeddings $\tau$ of $\mu$.
%
%th7.1 #&#
\begin{theorem}
\label{thmmain2}
%{\bf Do I not need $\bar{\Lambda}_\beta< \infty$ etc}
Suppose $g\dvtx  \R_+ \mapsto\R_+$ is a positive function and that $\mu$ is
centered.
\begin{longlist}[(ii)]
\item[(i)]
Suppose $g$ is increasing.
Then
\[
\inf_{\tau\in\sS(W,\mu)} \E\biggl[ \int_0^\tau
g(S_u) \,du \biggr] = \E\biggl[ \int_0^{\tau^{\mathrm{AY}}_\mu}
g(S_u) \,du \biggr]
\]
and
\[
\sup_{\tau\in\sS_{\mathrm{UI}}(W,\mu)} \E\biggl[ \int_0^\tau
g(S_u) \,du \biggr] = \E\biggl[ \int_0^{\tau^{\rmP }_\mu}
g(S_u) \,du \biggr].
\]

\item[(ii)]
Suppose $g$ is decreasing.
Then
\[
\inf_{\tau\in\sS(W,\mu)} \E\biggl[ \int_0^\tau
g(S_u) \,du \biggr] = \E\biggl[ \int_0^{\tau^{\rmP }_\mu}
g(S_u) \,du \biggr]
\]
and
\[
\sup_{\tau\in\sS_{\mathrm{UI}}(W,\mu)} \E\biggl[ \int_0^\tau
g(S_u) \,du \biggr] = \E\biggl[ \int_0^{\tau^{\mathrm{AY}}_\mu}
g(S_u) \,du \biggr].
\]
\end{longlist}
\end{theorem}
%
%re7.2 #&#
\begin{remark}
\label{remcounterintuitive}
As we remarked in the \hyperref[sec-intro]{Introduction}, at first sight this result is
counterintuitive. Given increasing $g$, the Az\'ema--Yor
stopping time maximizes $\E[ g(S_\tau) ]$ over $\tau\in
\sS_{\mathrm{UI}}(W,\mu)$, and it seems plausible that it might also maximize
$\E[ \int_0^\tau g(S_u)\,du ]$. In fact the exact opposite is true. The
explanation is that for the Az\'ema--Yor embedding there is
co-monotonicity\setcounter{footnote}{1}\footnote{A pair of random variables $X$ and $Y$ is
co-monotonic if $\Prob(X \leq x, Y\leq
y) = \min\{ \Prob(X \leq x), \Prob(Y \leq y) \}$ for all $x$ and $y$.}
between $S_\tau$
and
$W_\tau$, and conditional on
$S_\tau\geq s$, the stopping time occurs quite soon [and
certainly before $W$ drops below $\beta(s)$], whereas for the
Perkins embedding, conditional on
$S_\tau\geq s$, there are paths which will only be stopped when $W$
goes
below $\alpha^+(s)$. Thus,
for increasing $g$ when we wish to maximize the time (before $\tau$)
for which $S$ is large, this is best achieved by the Perkins embedding:
although relatively few paths will have large $S$ (most will have
already been
stopped) those with a large maximum will spend a long time after
first hitting $s$ before being
stopped.
\end{remark}
%
%ex7.3 #&#
\begin{example} Recall Example~\ref{expowerunif2}.
Suppose $\mu= U[-1,1]$ and $g(s) = s^{-c}$.
Then, for $c<0$, $(1 - c)^{-1}(2 - c)^{-1}(3 - c)^{-1} \leq\E[
\int_0^\tau S_u^{-c} \,du ] \leq\break (2 - c)^{-1}(3/2 - c)^{-1}$.

For $0<c<1$, $(2 - c)^{-1}(3/2 - c)^{-1} \leq\E[
\int_0^\tau S_u^{-c} \,du ] \leq(1 - c)^{-1}(2 - c)^{-1}(3 - c)^{-1}$,
for $1 \leq c < 3/2$, $(2 - c)^{-1}(3/2 - c)^{-1} \leq\E[
\int_0^\tau S_u^{-c} \,du ] \leq\infty$ and for $c \geq3/2$,
$\E[\int_0^\tau S_u^{-c} \,du ]=\infty$ for all embeddings $\tau$.

Note that for $c=0$, $\E[\tau]$ is independent of $\tau$ and equal to
the variance of~$\mu$.\vadjust{\goodbreak}
\end{example}
%
%ex7.4 #&#
\begin{example}
Recall the calculations from Example
\ref{excalcpareto}. Let the target law $\mu$ with support
$[-1,\infty)$
satisfy $\mu(dx)=\frac{2}{(x+2)^3}\,dx$. Let $g(s)=\frac{1}{c+s}$ for
$c>0$ which is decreasing in $s$.

The Az\'{e}ma--Yor upper bound can be calculated explicitly to be
\begin{eqnarray*}
B^{\mathrm{AY}}(c) &=& \int_{-1}^\infty
\frac{(b(w)-w)^2}{b(w)+c} \frac{2}{(w+2)^3} \,dw
\\
%&=&\int_0^\infty\frac{(s-\beta(s))^2}{s+c} \Prob(S_{\tau_\mu^{\mathrm{AY}}} \in
%ds) \\
%&=& \int_0^1 \frac{(s-(s/2-1))^2}{s+c} \,\frac{ds}{(1+s/2)^3} \\
&=& \frac{2(\log(c)-\log(2))}{c-2}.
\end{eqnarray*}

The expression for the Perkins lower bound is given by
\[
B^{\rmP }(c) = \int_{-1}^\infty
\frac{(a^+(w)-w)^2}{a^+(w)+c} \frac{2}{(w+2)^3} \,dw.
\]
%
%&=& \int_0^\infty\frac{(s-\alpha^+(s))^2}{s+c}
%&=& \int_0^\infty\frac{(s-\alpha^+(s))}{x+c} [1-(F(s)-F(\alpha^+(s)))]
%ds
The expression for $\alpha^+$ is too complicated for the expression
above to have an analytic representation. However, the values can be
computed
numerically for different $c$.
\end{example}

The rest of this section is devoted to a proof of
Theorem~\ref{thmmain2}. We split the proof into four separate parts.
\begin{pf*}{Proof of Theorem~\ref{thmmain2}\textup{(i)}: Lower bound}
Suppose first that $g$ is monotonic increasing and that we are
interested in minimizing the quantity $\E[ \int_0^\tau g(S_u) \,du ]$ over
embeddings $\tau$ of $\mu$ in $W$.
Note that it is sufficient to restrict attention to
$\sS_{\mathrm{UI}}(W,\mu)$: for nonminimal $\tau\in\sS(W,\mu)$ there exists
$\tilde{\tau} \leq\tau$ with $\tilde{\tau} \in\sS_{\mathrm{UI}}(W,\mu)$,
and then $\int_0^\tau g(S_u) \,du \geq\int_0^{\tilde{\tau}} g(S_u)
\,du$ for
each
$\omega\in\Omega$.

Suppose temporarily that $g$ is bounded and continuously differentiable.
Later we will relax this assumption.
Then $G(w,s) = (s-w)^2 g(s)$ satisfies G-MON$\downarrow$.

%Let $F_\mu$ denote the distribution function of $\mu$ and let $k_n
%= F^{-1}_\mu(1- 1/n)$.
For $\tau\in\sS_{\mathrm{UI}}(W,\mu)$ let $\sigma_n
= \tau\wedge H_{\pm n}$, let $\mu_n = \sL(W_{\sigma_n})$,
$\beta_n$ be the inverse barycenter of $\mu_n$ and finally let
$\tau^{\mathrm{AY}}_{\mu_n}$ be the Az\'{e}ma--Yor stopping rule associated with
the law $\mu_n$ so that $\tau^{\mathrm{AY}}_{\mu_n} = \tau_{\beta_n} =
\inf
\{ u\dvtx  W_u
\leq
\beta_n(S_u) \}$.
Then, by Proposition~\ref{propmuUbeta}, since $U_{\mu_n} \uparrow
U_{\mu}$,
%With the notation as above, $\mu_n$ converges weakly to $\mu$,
%$U_{\mu_n} \rightarrow U$ uniformly in $x$, $\beta_n(s) \rightarrow
$\tau_{\beta_n} \rightarrow\tau_\beta$ almost surely.

If a stopping time $\rho$ is such that $\rho\leq H_{\pm n}$, then
$\E[\rho]<\infty$ and for $u \leq\rho$, $(S_u - W_u)g(S_u)$ is
bounded. Then if $M_t = \int_0^t (S_u - W_u)g(S_u) \,dW_u$, we have that
$(M^\rho_t)_{t \geq0}$ is an $L^2$ bounded martingale for which
%
%e7.2 #&#
\begin{equation}
\label{Minfinity=0} \E\bigl[M^\rho_\infty\bigr] = \E\biggl[
\int_0^{\rho}(S_u -
W_u)g(S_u) \,du \biggr] = 0.
\end{equation}
It follows that
\begin{eqnarray*}
\E\biggl[ \int_0^{\sigma_n} g(S_u) \,du
\biggr] & = & \E\bigl[ (S_{\sigma_n} - W_{\sigma_n})^2
g(S_{\sigma_n}) \bigr]
\\
& \geq& \E\bigl[ (S_{\tau_{\beta_n}} - W_{\tau_{\beta_n}})^2
g(S_{\tau_{\beta_n}}) \bigr]
\\
& = & \E\biggl[ \int_0^{\tau_{\beta_n}}
g(S_u) \,du \biggr],
\end{eqnarray*}
where we have used (\ref{eqnitoG}) and (\ref{Minfinity=0}) twice
and
Theorem~\ref{thmmain1b}.
Then it follows from the Fatou lemma that
%
%e7.3 #&#
\begin{eqnarray}
\label{eqn51i} \lim\E\biggl[ \int_0^{\sigma_n}
g(S_u) \,du \biggr] &\geq&\lim\E\biggl[ \int_0^{\tau_{\beta_n}}
g(S_u) \,du \biggr] \nonumber\\[-8pt]\\[-8pt]
&\geq&\E\biggl[ \linf\int_0^{\tau_{\beta_n}}
g(S_u) \,du \biggr]\nonumber
\end{eqnarray}
and by monotone convergence and the fact that
$\tau_{\beta_n}
\rightarrow\tau_\beta$ almost surely,
\[
\E\biggl[ \int_0^{\tau}
g(S_u) \,du \biggr]
\geq\E\biggl[\int_0^{\tau_\beta} g(S_u) \,du \biggr]
\]
as required.

Finally we remove the temporary assumptions on $g$. Given $g$
is monotonic increasing we can find an increasing sequence of
bounded, continuously differentiable (increasing) functions $g_m$ which
approximate $g$ from below. Then, by monotone convergence
\begin{eqnarray*}
\E\biggl[ \int_0^{\tau} g(S_u) \,du
\biggr] &=& \lim_m \E\biggl[ \int_0^{\tau}
g_m (S_u) \,du \biggr] \geq\lim_m \E\biggl[
\int_0^{\tau_\beta} g_m (S_u)
\,du \biggr] \\
&=& \E\biggl[ \int_0^{\tau_\beta}
g(S_u) \,du \biggr].
\end{eqnarray*}

Note that this same argument will apply in all four parts of
Theorem~\ref{thmmain2}, and henceforth without loss of generality we
will assume that $g$ is continuously differentiable and bounded by
$\bar{g}$.
\noqed\end{pf*}
\begin{pf*}{Proof of Theorem~\ref{thmmain2}\textup{(ii)}: Lower bound}
\textit{Case} 1: There exists an open interval $I \subseteq
[-1,1]$ containing 0
with $\mu(I)=0$.

Given $\tau\in\sS(W,\mu)$, let $\sigma_m = \tau\wedge H_{\pm m}$. Let
$\mu_m = \sL(W_{\sigma_m})$. Write $\tau^{\rmP }_m$ for the Perkins embedding
of $\mu_m$. Note that $\mu_m \Rightarrow\mu$, $U_{\mu_m}(0)
\rightarrow
U_{\mu}(0)$ and $\mu_{m}(I)=0$. Then, $\tau^{\rmP }_m = \tau_{\alpha_m}$ and
by Proposition~\ref{propmuUalpha}(a), $\tau_{\alpha_m} \rightarrow
\tau_{\alpha}$ almost surely. Then
exactly as in (\ref{eqn51i}), but now using Theorem~\ref{thmmain1a}
to give that the lower bound is attained by the Perkins embedding, we
conclude that $\E[ \int_0^{\tau} g(S_u) \,du ] \geq
\E[ \int_0^{\tau^{\rmP }_\mu} g(S_u) \,du ]$.

\textit{Case} 2: General $\mu$.
Given any subsequence, by Proposition~\ref{propmuUalpha}(b) we may
take a
further subsequence down which
$\tau^{\rmP }_m \rightarrow\tau^{\rmP }$ almost surely.
Then down this subsequence the result holds, as in case 1. Since the
first subsequence was arbitrary we are done.
\noqed\end{pf*}
\begin{pf*}{Proof of Theorem~\ref{thmmain2}\textup{(ii)}: Upper bound}
Now consider the upper bound in Theorem~\ref{thmmain2}(ii). Rather
than attempting to find a dominating random variable which will
allow us to use the reverse Fatou lemma in place of the Fatou lemma
above
we will use a slightly different approach based on defining a
sequence of intermediate stopping times.

Let $\tau$ be any element of $\sS_{\mathrm{UI}}(W,\mu)$. Suppose $g$ is
bounded, continuously differentiable and
monotonic decreasing, and that $\mu$ has support in a bounded interval
$[\check{x},\hat{x}]$.
Then, as above, $\E[ \int_0^{\tau} g (S_u) \,du ] =
\E[G(W_\tau,S_\tau)]$.
Moreover, we can conclude
from Theorem~\ref{thmmain1a} that
\[
\sup_{\tau\in\sS_{\mathrm{UI}}(W,\mu)} \E\biggl[ \int_0^{\tau}
g(S_u) \,du \biggr] = \E\biggl[ \int_0^{\tau_\beta}
g(S_u) \,du \biggr].
\]

It remains to remove the assumptions on $\mu$.

Given $\varepsilon$, let $U_\varepsilon(x) = \max\{ U_{\mu}(x)-\varepsilon, |x|
\}$, and let $\check{x}_\varepsilon$ and $\hat{x}_\varepsilon$ be the
left and
right-hand
endpoints of the interval $I_{\varepsilon} =\{ x\dvtx  U_{\varepsilon}(x)>|x| \}$.

Let $\sigma_\varepsilon= \tau\wedge\inf\{ u\dvtx  W_u \notin I_{\varepsilon}
\}$. Let $\bar{\mu}_\varepsilon$ be the law of $W_{\sigma_\varepsilon}$, and
let $\bar{U}_\varepsilon$ be the associated potential. Then
$\bar{U}_\varepsilon= U_\varepsilon$ on $I_\varepsilon^c$ and
$U_\varepsilon\leq\bar{U}_\varepsilon\leq U_\mu$.

%f5 #&#
\begin{figure}

\includegraphics{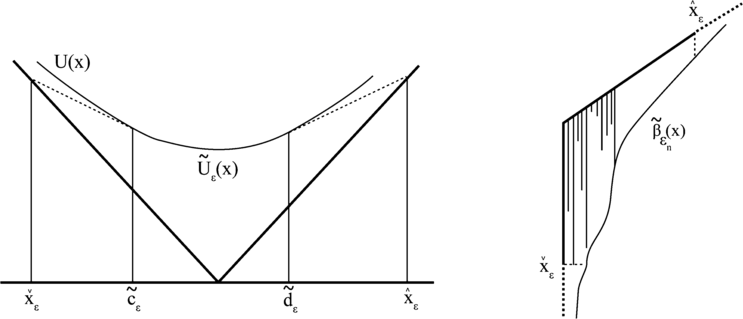}

\caption{The potentials $\tilde{U}_\varepsilon$ increase monotonically as
$\varepsilon$ decreases. Moreover, over a range of $x$, depending on
$\varepsilon_n$, we have
$\tilde{\beta}_{\varepsilon_n}(x) \equiv\beta(x)$, and hence,
the inverse barycentre functions converge monotonically.}\label
{figfinaldouble}
\end{figure}

Now let $\tilde{U}_\varepsilon$ be the largest convex function such that
$\tilde{U}_\varepsilon(x) = |x|$ on $I^c_{\varepsilon}$ and
$\tilde{U}_\varepsilon
\leq U_\mu$. It follows that $\tilde{U}_\varepsilon$ is actually equal to
$U$
on an interval $\tilde{I}_\varepsilon= [\tilde{c}_\varepsilon,
\tilde{d}_\varepsilon]$. If $\varepsilon$ is small enough, then $0 \in
\tilde{I}_{\varepsilon}$.
See Figure~\ref{figfinaldouble}.
Further, $U_\varepsilon\leq\bar{U}_\varepsilon\leq\tilde{U}_\varepsilon
\leq
U$
and in terms of the
associated measures
$\mu_\varepsilon\leq_{\mathrm{cx}} \bar{\mu}_\varepsilon\leq_{\mathrm{cx}}
\tilde{\mu}_\varepsilon
\leq_{\mathrm{cx}} \mu$, where $\tilde{\mu}_\varepsilon$ is such that
$U_{\tilde{\mu}_\varepsilon} = \tilde{U}_\varepsilon$, and we
recall that
$\leq_{\mathrm{cx}}$ denotes ``less than or equal to in convex order.''
Then, by a theorem of Strassen~\cite{Strassen67} (or for a more
explicit construction in our context, Chacon and
Walsh~\cite{ChaconWalsh76}), given $\sigma_\varepsilon$
there exists a
stopping time $\tilde{\sigma}_\varepsilon$ such
that $\sigma_\varepsilon\leq\tilde{\sigma}_\varepsilon$ almost surely, and
$\tilde{\mu}_\varepsilon=
\sL(W_{\tilde{\sigma}_\varepsilon})$.

Now consider a sequence $\varepsilon_n$ decreasing to zero.
Let $\tilde{\beta}_{\varepsilon_n}$ be the inverse barycentre associated
with
$\tilde{\mu}_{\varepsilon_n}$, and let $\tilde{\tau}_n$ be the Az\'{e}ma--Yor
stopping time associated with~$\tilde{\beta}_{\varepsilon_n}$. The
introduction of
the stopping times $\tilde{\sigma}_{\varepsilon_n}$ gives extra
structure which
means
that not only do the barycenters converge (as in
Proposition~\ref{propmuUbeta}), but also
that they
converge monotonically.
\noqed\end{pf*}
%
%le7.5 #&#
\begin{lemma}
$\tilde{\beta}_n \downarrow\beta$ and $\tilde{\tau}_n \uparrow
\tau_\beta$ almost surely.
\end{lemma}
\begin{pf}
Write $\check{x}_n$ (resp., $\hat{x}_n, c_n, d_n$) for
$\check{x}_{\varepsilon_n}$
(resp.,
$\hat{x}_{\varepsilon_n}, c_{\varepsilon_n}, d_{\varepsilon_n}$).

Then, for $s \leq b({c}_n)$, $\tilde{\beta}_n(s)= \check{x}_n \geq
\beta(s)$,
for
$b(\tilde{c}_n) < s < \hat{x}_n$, $\tilde{\beta}_n(s)=
{\beta}(s)$ and for
$s \geq\hat{x}_n$, $\tilde{\beta}_n(s)= s \geq\beta(s)$.

Monotonicity in $n$ of $\tilde{\tau}_n$ follows immediately.
\end{pf}

It follows from the results for bounded target distributions that
\begin{eqnarray*}
\E\biggl[ \int_0^{\sigma_n}
g(S_u) \,du \biggr] &\leq&\E\biggl[ \int_0^{\tilde{\sigma}_n}
g(S_u) \,du \biggr] = \E\bigl[ G(W_{\tilde{\sigma}_n},S_{\tilde{\sigma}_n})
\bigr] \leq\E\bigl[ G(W_{\tilde{\tau}_n},S_{\tilde{\tau}_n}) \bigr] \\
&=& \E\biggl[ \int_0^{\tilde{\tau}_{n}} g(S_u) \,du \biggr].
\end{eqnarray*}
We have that the integral inside the first expectation converges
monotonically
to $\int_0^{\tau} g(S_u) \,du$, whereas the integral inside the
final expression converges monotonically to
$\int_0^{\tau_\beta} g(S_u) \,du$.
Hence $\E[ \int_0^{\tau} g(S_u) \,du ] \leq\E[ \int_0^{\tau_\beta}
g(S_u) \,du ]$ as required.\vspace*{-2pt}
\begin{pf*}{Proof of Theorem~\ref{thmmain2}\textup{(i)}: Upper bound}
The final element of Theorem~\ref{thmmain2} is the upper bound in
the case of monotonically increasing $g$.
Recall that we suppose that $g$ is continuously differentiable,
and bounded by
$\bar{g}$.

If $\mu$ has bounded support, then Theorem~\ref{thmmain1b} applies
directly, so we assume that the support of $\mu$ is unbounded.

If $\mu\notin L^2$, then for each $\tau\in\sS(W,\mu)$ we have $\E
[\tau]
=
\infty$ and
using the fact that $\E[{H_{\varepsilon} \wedge\tau^{\rmP }_\mu}] \leq\E
[H_{\varepsilon}
\wedge H_{\alpha^+(\varepsilon)}]<\infty$, we have that
$\E[ \int_0^{\tau^{\rmP }_\mu} g(S_u) \,du ] \geq g(\varepsilon) \E[
\int_{H_{\varepsilon} \wedge\tau^{\rmP }_\mu}^{\tau^{\rmP }_\mu} \,du ] = \infty
$, and
there is
nothing to prove.

So suppose $\mu\in L^2$. Then the area between the curves $U_{\mu}(x)$
and $|x|$ is finite.

Let $U_\varepsilon(x) = \max\{ U_{\mu}(x)-\varepsilon, |x|
\}$ and related quantities be defined as above.

This time,
since $\tilde{U}_\varepsilon\equiv U_\mu$ on $\tilde{I}_\varepsilon$ we have
that $\alpha_{\tilde{\mu}_\varepsilon} = \alpha_\mu$ on some sub-interval
$\acute{I}_\varepsilon\subseteq\tilde{I}_\varepsilon$ of the form
$\acute{I}_\varepsilon= [\acute{c}_\varepsilon,\acute{d}_\varepsilon]$,
and as
$\varepsilon
\downarrow0$, $\acute{I}_\varepsilon$ increases to the support of~$\mu$.

Now
\[
\E\biggl[ \int_0^{\tau} g(S_u) \,du
\biggr] = \lim_{\varepsilon
\downarrow
0} \E\biggl[ \int_0^{\sigma_\varepsilon}
g(S_u) \,du \biggr]
\]
and
\[
\E\biggl[ \int_0^{\sigma_\varepsilon} g(S_u) \,du
\biggr] \leq\E\biggl[ \int_0^{\tilde{\sigma}_\varepsilon}
g(S_u) \,du \biggr] \leq\E\biggl[ \int_0^{\tau^{\rmP }(\tilde{\mu}_\varepsilon)}
g(S_u) \,du \biggr].
\]

But
\begin{eqnarray*}
\E\biggl[ \int_0^{\tau^{\rmP }(\tilde{\mu}_\varepsilon)} g(S_u) \,du
\biggr] &=& \E\biggl[ \int_0^{\tau^{\rmP }(\tilde{\mu}_\varepsilon) \wedge
H_{\acute{c}_\varepsilon}
\wedge H_{\acute{d}_\varepsilon}}
g(S_u) \,du \biggr]\\[-2pt]
&&{}+ \E\biggl[ \int_{\tau^{\rmP }(\tilde{\mu}_\varepsilon) \wedge
H_{\acute{c}_\varepsilon}
\wedge H_{\acute{d}_\varepsilon}}^{\tau^{\rmP }(\tilde{\mu}_\varepsilon)}
g(S_u) \,du \biggr].
\end{eqnarray*}
Since $\alpha_{\tilde{\mu}_\varepsilon}= \alpha_\mu$ on
$\acute{I}_\varepsilon$
and we have that $\tau^{\rmP }(\tilde{\mu}_\varepsilon) \wedge
H_{\acute{c}_\varepsilon}\wedge H_{\acute{d}_\varepsilon}$ is monotonically
increasing as $\varepsilon\downarrow0$ and hence
the first term on the right-hand side converges to
$\E[ \int_0^{\tau^{\rmP }({\mu})} g(S_u) \,du ]$.
Meanwhile, the second term is bounded by
$\bar{g} \E[ \tau^{\rmP }(\tilde{\mu}_\varepsilon) - \tau^{\rmP }(\tilde{\mu
}_\varepsilon)
\wedge H_{\acute{c}_\varepsilon} \wedge H_{\acute{d}_\varepsilon}]$. This last
quantity
is at
most
$\bar{g}$ multiplied by the area between the potentials
$U_\mu$ and $U_{\acute{\mu}_\varepsilon}$ where $\acute{\mu
}_\varepsilon=
\sL(W_{\tau^{\rmP }(\tilde{\mu}_\varepsilon) \wedge H_{\acute{c}} \wedge
H_{\acute{d}}})$.
However, as $\varepsilon$ tends to zero this area tends to zero.
Hence $\E[ \int_0^{\tau} g(S_u) \,du ] \leq
\E[ \int_0^{\tau^{\rmP }(\mu)} g(S_u) \,du ]$.
\end{pf*}

%s8 #&#
\section{An application and extensions}
%s8.1 #&#
\subsection{Variance swap on the sum of squared returns}
\label{ssecvswap}

We now return to the question which originally motivated this paper
which was to
find model-independent bounds for variance swaps given the terminal law
of the
underlying asset price process or equivalently, call prices with expiry
$T$ for all
strikes. Using the results developed in this article we will show how
to bound the
idealized variance swap based on squared returns, introduced in Section~\ref{sec-squaredreturns}. The results in this article motivated
further work on
model-independent bounds and hedging strategies for variance swaps in a general
setting; see Hobson and Klimmek~\cite{HobsonKlimmek10}.
%
%The relationship between
%variance swap bounds and the Skorokhod embedding problem solved in
%Section

As in Section~\ref{sec-squaredreturns}, let $X=(X_t)_{0 \leq t \leq T}$
be a square-integrable martingale started at $X_0=x_0$ with $X_T \sim
\mu$,
where $\mu$ is centered at $x_0$ and supported on $\R_+$. Recall from
(\ref{eqidealised}) the
definition
for the payoff of an idealized variance swap
$V_T=V((X_s)_{0 \leq s \leq T})=\int_0^T {(X_{t-})^{-2}} \,{d[X,X]_t}$.
By (\ref{eqvarlower}) and (\ref{eqvarupper}) we have
\[
\inf_{\tau\in S_{\mathrm{UI}}(B,\mu)} \E\biggl[\int_0^\tau
\frac
{du}{(S^B_u)^2} \biggr] \leq\E[V_T] \leq\sup_{\tau\in S_{\mathrm{UI}}(B,\mu)} \E
\biggl[\int_0^\tau\frac{du}{(I^B_u)^2} \biggr].
\]
Let $\tilde{\mu}$ be the measure $\mu$ reflected around $0$, so that
$\tilde\mu$
is a
measure on $\R_-$, and observe that
\[
\sup_{\tau\in S_{\mathrm{UI}}(B,\mu)} \E\biggl[\int_0^\tau
\frac
{du}{(I^B_u)^2} \biggr] = \sup_{\tau\in S_{\mathrm{UI}}(\tilde{B},\tilde\mu)} \E
\biggl[\int
_0^\tau\frac{du}{(S^{\tilde{B}}_u)^2} \biggr],
\]
where $\tilde{B}$ is a Brownian motion started at $-x_0$, with maximum process
$S^{\tilde{B}}$.
Now we apply Theorem~\ref{thmmain2} to see that
\[
\E\biggl[\int_0^{\tau_\mu^{\rmP }} \frac{du}{(S^{B}_u)^2}
\biggr] \leq\E[V_T] \leq\E\biggl[\int_0^{\tau_{\tilde{\mu}}^{\rmP }}
\frac{du}{(S^{\tilde
{B}}_u)^2} \biggr].
\]
Note that the Perkins embedding for $\tau_{\tilde\mu}$ is determined
by the
monotonic
functions $\alpha_{\tilde\mu}^\pm$ where
$\alpha_{\tilde\mu}^\pm(x)=-\alpha^\mp_{\mu}(-x)$.
%Moreover, it can be shown that there is a martingale for which the
%bound is attained.
%
%ex8.1 #&#
\begin{example}
Suppose that $X_0=1$ and $\mu=U[0,2]$.
Shifting the quantities calculated in Example~\ref{excalcuniform} to
allow for the starting value $X_0=1$ it is clear that $\alpha_\mu^+\dvtx
[1,2] \rightarrow[0,1]$ is defined
$\alpha_\mu^+(s)=s-2\sqrt{s-1}$ and $\alpha_{\mu}^-\dvtx  [0,1]
\rightarrow[1,2]$ is defined $\alpha_{\mu}^-(i)=i+\sqrt{1-i}$.
Hence the lower bound can be calculated,
\[
\E\biggl[\int_0^{\tau^{\rmP }_\mu} \frac{du}{S_u^2}
\biggr] = \E\biggl[ \biggl( 1 - \frac{B_{\tau^{\rmP }_\mu}}{S_{\tau^{\rmP }_\mu}}
\biggr)^2 \biggr]
= \int_0^1 \biggl( 1 - \frac{x}{a_\mu^+(x)}
\biggr)^2 \,\frac{dx}{2} %= \int_1^2 \frac{(s-\alpha^+(s))^2}{s^2} \,d
= \frac{\pi}{2}-2
\log{2}.
\]

For the upper bound, first considering $g_{\varepsilon}(s) = s^{-2}
\wedge
\varepsilon^{-2}$ and then letting $\varepsilon\downarrow0$,
\[
\E\biggl[\int_0^{\tau^{\rmP }_{\tilde\mu}} \frac{du}{\tilde
{S}_u^2}
\biggr] = \E\biggl[ \biggl( 1 - \frac{\tilde{B}_{\tau^{\rmP }_{\tilde\mu
}}}{\tilde{S}_{\tau^{\rmP }_{\tilde
{\mu}}}} \biggr)^2 \biggr]
= \int_0^1 \biggl( 1 - \frac{x}{a_\mu^-(x)}
\biggr)^2 \,\frac{dx}{2} %= \int_1^2 \frac{(s-\alpha^+(s))^2}{s^2} \,d
= \infty.
\]
\end{example}

%s8.2 #&#
\subsection{Extension to diffusions}

Suppose that $(X_t)_{t \geq0}$ is a time-homogeneous
diffusion on $I \subseteq\R$.
More specifically, let $\sigma\dvtx I \rightarrow(0,\infty)$ and $b\dvtx  I
\rightarrow\R$ be Lipschitz functions, and define $(X_t)_{t \geq0}$ to
be the solution to
\[
dX_t= \sigma(X_t)\,dB_t +
b(X_t) \,dt,\qquad  X_0=x_0,
\]
where $(B_t)_{t \geq0}$ is a Brownian motion.

Let $s\dvtx  I \rightarrow\R$ be the strictly increasing and $C^2$ scale
function of $X$,
\[
s(x_0)=0,\qquad s'(x)=\exp\biggl(-\int
_0^x 2 \frac{b(u)}{\sigma(u)^{2}} \,du \biggr)
\]
and let $h=s^{-1}$.

Consider the problem of maximizing (or minimizing)
$\E[F(X_\tau,S^X_\tau)]$ over minimal embeddings $\tau$ of $\mu$.
Since $M_t =
s(X_t)$ is a local martingale it follows that it can be represented
as $M_t=W_{A(t)}$, for some (continuous) time-change $t \rightarrow
A(t)$. Define the
measure $\nu$ by $\nu(G)=\mu(s^{-1}(G))$ for Borel sets $G \subseteq
s(I)$.
Notice that $\sigma$ is a minimal embedding of $\nu$ in $W$ if and only
if $\tau=A^{-1}(\sigma)$ is a minimal embedding of $\nu$ in $M$ and
hence a minimal embedding of $\mu$ in~$X$.

%We refer to Cox and Hobson
%property
%is carried through a transformation by a scale function when its range
%is a strict subset of $\R$.

Define the function $\Fhat$ by $\Fhat(w,s)=F(h(w),h(s))$. Then
%
%e8.1 #&#
\begin{equation}
\label{eqscaledF} F\bigl(X_\tau,S^X_\tau\bigr)
= \F\bigl(h(W_{A_\tau}),h(S_{A_\tau})\bigr) = \Fhat(W_{A_\tau},S_{A_\tau}).
\end{equation}
%
%le8.2 #&#
\begin{lemma} \label{ladditions-FMON}
Suppose $F$ satisfies F-MON$ \uparrow$. Then $\hat{F}$ satisfies
$\hat{F}$-MON$ \uparrow$ if $F_s < 0$ and $h$ is concave or if
$F_s > 0$ and $h$ is convex.

Similarly, suppose $F$ satisfies F-MON$ \downarrow$. Then $\hat{F}$
satisfies $\hat{F}$-MON$ \downarrow$ if $F_s < 0$ and $h$ is convex or if
$F_s>0$ and $h$ is concave.
\end{lemma}
\begin{pf}
The result follows from the expression
%
%e8.2 #&#
\begin{equation}
\label{eqscaleFMON} \frac{\Fhat_s(x,s)}{s-x} = \frac{h'(s)
F_s(h(x),h(s))}{h(s)-h(x)} \frac{h(s)-h(x)}{s-x}.
\end{equation}
\upqed\end{pf}

Note that $h$ is convex (concave) when $s$ is concave (convex), and
since $2s''(x)/s'(x)=-\sigma(x)^2/b(x)$, the scale function is
concave if $b(x) >0$ for all~$x$.
%
%pr8.3 #&#
\begin{proposition} \label{pdiffusioncase}
Suppose $\nu= \mu\circ h$ is centered about zero, and suppose $b>0$.
Suppose $F$ satisfies
F-MON$\uparrow$ and is increasing in $s$. Then
\begin{eqnarray*}
\sup_{\tau\in S_{\mathrm{UI}}(X,\mu)} \E\bigl[F\bigl(X_\tau,S^X_\tau
\bigr)\bigr] &=&\E\bigl[\Fhat(W_{\tau_\nu^{\mathrm{AY}}},S_{\tau_\nu^{\mathrm{AY}}})\bigr
],
\\
\inf_{\tau\in S_{\mathrm{UI}}(X,\mu)} \E\bigl[F\bigl(X_\tau,S^X_\tau
\bigr)\bigr] &=&\E\bigl[\Fhat(W_{\tau_\nu^{\rmP }},S_{\tau_\nu^{\rmP }})\bigr].
\end{eqnarray*}
%
%The result is similar for the other cases listed in Lemma
\end{proposition}
%
%re8.4 #&#
\begin{remark}
Whilst necessary to apply the results of the Brownian setting, the
assumption that $\nu\equiv\mu\circ h$ is centered is not as innocuous
as might first appear, and in the setting of a transient diffusion it is
natural to wish to consider embeddings for target laws which, after
transformation by the scale function, are not centered. For example, let
$X$ be a three-dimensional Bessel process, started at one.
Then $s(x) = -1/x + 1$ and $h(m) = 1/(1-m)$.
Now let $\mu$
be any probability measure on $\R^+$ with $\int_{\R^+} x^{-1} \mu(dx)
\leq1$. Then, there exists
a minimal embedding of $\mu$ in $X$, but only if $\int_{\R^+} x^{-1}
\mu(dx)=1$ does this embedding correspond to a uniformly integrable
embedding of $M \equiv1 - X^{-1}$.

See Cox and Hobson~\cite{CoxHobson06} (and the references therein) for
a further discussion of this issue, and of the construction of
embeddings in Brownian motion of noncentered target laws.
\end{remark}

\section*{Acknowledgments}
Both authors thank a pair of referees for their detailed comments on an
earlier version of this article.

%suskaldyti doi

% imsref loaded by lrinkeviciute, 2013-01-29 11:20:21

\printaddresses


\begin{thebibliography}{24}
% BibTex style file: ims.bst, 2013-01-28
% Default style options (sort=0,type=number).
% Used options (sort=1,type=number).

%b1 ###
\bibitem{AzemaYor79b}
\begin{bincollection}[mr]
\bauthor{\bsnm{Az{\'e}ma},~\bfnm{Jacques}\binits{J.}} \AND
  \bauthor{\bsnm{Yor},~\bfnm{Marc}\binits{M.}}
(\byear{1979}).
\btitle{Le probl\`eme de {S}korokhod: Compl\'ements \`a ``{U}ne solution simple
  au probl\`eme de {S}korokhod.''}
In \bbooktitle{S\'eminaire de {P}robabilit\'es, {XIII} ({U}niv. {S}trasbourg,
  {S}trasbourg, 1977/78)}.
\bseries{Lecture Notes in Math.}
\bvolume{721}
\bpages{625--633}.
\bpublisher{Springer}, \blocation{Berlin}.
\bid{mr={0544832}}
\bptok{imsref}%
\end{bincollection}
\endbibitem

%b2 ###
\bibitem{AzemaYor79a}
\begin{bincollection}[mr]
\bauthor{\bsnm{Az{\'e}ma},~\bfnm{Jacques}\binits{J.}} \AND
  \bauthor{\bsnm{Yor},~\bfnm{Marc}\binits{M.}}
(\byear{1979}).
\btitle{Une solution simple au probl\`eme de {S}korokhod}.
In \bbooktitle{S\'eminaire de {P}robabilit\'es, {XIII} ({U}niv. {S}trasbourg,
  {S}trasbourg, 1977/78)}.
\bseries{Lecture Notes in Math.}
\bvolume{721}
\bpages{90--115}.
\bpublisher{Springer}, \blocation{Berlin}.
\bid{mr={0544782}}
\bptok{imsref}%
\end{bincollection}
\endbibitem

%b3 ###
\bibitem{BreedenLitzenberger78}
\begin{barticle}[auto:STB|2013/01/23|16:20:06]
\bauthor{\bsnm{Breeden},~\bfnm{D.~T.}\binits{D.~T.}} \AND
  \bauthor{\bsnm{Litzenberger},~\bfnm{R.~H.}\binits{R.~H.}}
(\byear{1978}).
\btitle{Prices of state-contingent claims implicit in option prices}.
\bjournal{J. Bus.}
\bvolume{51}
\bpages{621--651}.
\bptok{imsref}%
\end{barticle}
\endbibitem

%b4 ###
\bibitem{Chacon77}
\begin{barticle}[mr]
\bauthor{\bsnm{Chacon},~\bfnm{R.~V.}\binits{R.~V.}}
(\byear{1977}).
\btitle{Potential processes}.
\bjournal{Trans. Amer. Math. Soc.}
\bvolume{226}
\bpages{39--58}.
\bid{issn={0002-9947}, mr={0501374}}
\bptok{imsref}%
\end{barticle}
\endbibitem

%b5 ###
\bibitem{ChaconWalsh76}
\begin{bincollection}[mr]
\bauthor{\bsnm{Chacon},~\bfnm{R.~V.}\binits{R.~V.}} \AND
  \bauthor{\bsnm{Walsh},~\bfnm{J.~B.}\binits{J.~B.}}
(\byear{1976}).
\btitle{One-dimensional potential embedding}.
In \bbooktitle{S\'eminaire de {P}robabilit\'es, {X} ({P}r\`emiere Partie,
  {U}niv. {S}trasbourg, {S}trasbourg, Ann\'ee Universitaire
  1974/1975)}.
\bseries{Lecture Notes in Math.}
\bvolume{511}
\bpages{19--23}.
\bpublisher{Springer}, \blocation{Berlin}.
\bid{mr={0445598}}
\bptok{imsref}%
\end{bincollection}
\endbibitem

%b6 ###
\bibitem{CoxHobson06}
\begin{barticle}[mr]
\bauthor{\bsnm{Cox},~\bfnm{A.~M.~G.}\binits{A.~M.~G.}} \AND
  \bauthor{\bsnm{Hobson},~\bfnm{D.~G.}\binits{D.~G.}}
(\byear{2006}).
\btitle{Skorokhod embeddings, minimality and non-centred target distributions}.
\bjournal{Probab. Theory Related Fields}
\bvolume{135}
\bpages{395--414}.
\bid{doi={10.1007/s00440-005-0467-y}, issn={0178-8051}, mr={2240692}}
\bptok{imsref}%
\end{barticle}
\endbibitem

%b7 ###
\bibitem{DemeterfiDerman99}
\begin{barticle}[auto:STB|2013/01/23|16:20:06]
\bauthor{\bsnm{Demeterfi},~\bfnm{K.}\binits{K.}},
  \bauthor{\bsnm{Derman},~\bfnm{E.}\binits{E.}},
  \bauthor{\bsnm{Kamal},~\bfnm{M.}\binits{M.}} \AND
  \bauthor{\bsnm{Zou},~\bfnm{J.}\binits{J.}}
(\byear{1999}).
\btitle{A guide to volatility and variance swaps}.
\bjournal{The Journal of Derivatives}
\bvolume{6}
\bpages{9--32}.
\bptok{imsref}%
\end{barticle}
\endbibitem

%b8 ###
\bibitem{Durrett91}
\begin{bbook}[mr]
\bauthor{\bsnm{Durrett},~\bfnm{Richard}\binits{R.}}
(\byear{1991}).
\btitle{Probability: Theory and Examples}.
\bpublisher{Wadsworth \& Brooks/Cole Advanced Books \& Software},
  \blocation{Pacific Grove, CA}.
\bid{mr={1068527}}
\bptok{imsref}%
\end{bbook}
\endbibitem

%b9 ###
\bibitem{HirschProfetaRoynetteYor10}
\begin{bincollection}[mr]
\bauthor{\bsnm{Hirsch},~\bfnm{Francis}\binits{F.}},
  \bauthor{\bsnm{Profeta},~\bfnm{Christophe}\binits{C.}},
  \bauthor{\bsnm{Roynette},~\bfnm{Bernard}\binits{B.}} \AND
  \bauthor{\bsnm{Yor},~\bfnm{Marc}\binits{M.}}
(\byear{2011}).
\btitle{Constructing self-similar martingales via two {S}korokhod embeddings}.
In \bbooktitle{S\'eminaire de {P}robabilit\'es {XLIII}}.
\bseries{Lecture Notes in Math.}
\bvolume{2006}
\bpages{451--503}.
\bpublisher{Springer}, \blocation{Berlin}.
\bid{doi={10.1007/978-3-642-15217-7_21}, mr={2790387}}
\bptok{imsref}%
\end{bincollection}
\endbibitem

%b10 ###
\bibitem{Hobson10}
\begin{bincollection}[mr]
\bauthor{\bsnm{Hobson},~\bfnm{David}\binits{D.}}
(\byear{2011}).
\btitle{The {S}korokhod embedding problem and model-independent bounds for
  option prices}.
In \bbooktitle{Paris-{P}rinceton {L}ectures on {M}athematical {F}inance 2010}.
\bseries{Lecture Notes in Math.}
\bvolume{2003}
\bpages{267--318}.
\bpublisher{Springer}, \blocation{Berlin}.
\bid{doi={10.1007/978-3-642-14660-2_4}, mr={2762363}}
\bptnote{check year}%
\bptok{imsref}%
\end{bincollection}
\endbibitem

%b11 ###
\bibitem{HobsonKlimmek10}
\begin{bmisc}[auto:STB|2013/01/23|16:20:06]
\bauthor{\bsnm{Hobson},~\bfnm{D.~G.}\binits{D.~G.}} \AND
  \bauthor{\bsnm{Klimmek},~\bfnm{M.}\binits{M.}}
(\byear{2012}).
\bhowpublished{Model-independent hedging strategies for variance swaps.
\textit{Finance Stoch.} DOI:\doiurl{10.1007/s00780-012-0190-3}}.
\bptok{imsref}%
\end{bmisc}
\endbibitem

%b12 ###
\bibitem{HobsonPedersen02}
\begin{barticle}[mr]
\bauthor{\bsnm{Hobson},~\bfnm{David~G.}\binits{D.~G.}} \AND
  \bauthor{\bsnm{Pedersen},~\bfnm{J.~L.}\binits{J.~L.}}
(\byear{2002}).
\btitle{The minimum maximum of a continuous martingale with given initial and
  terminal laws}.
\bjournal{Ann. Probab.}
\bvolume{30}
\bpages{978--999}.
\bid{doi={10.1214/aop/1023481014}, issn={0091-1798}, mr={1906424}}
\bptok{imsref}%
\end{barticle}
\endbibitem

%b13 ###
\bibitem{KertzRosler90}
\begin{barticle}[mr]
\bauthor{\bsnm{Kertz},~\bfnm{Robert~P.}\binits{R.~P.}} \AND
  \bauthor{\bsnm{R{\"o}sler},~\bfnm{Uwe}\binits{U.}}
(\byear{1990}).
\btitle{Martingales with given maxima and terminal distributions}.
\bjournal{Israel J. Math.}
\bvolume{69}
\bpages{173--192}.
\bid{doi={10.1007/BF02937303}, issn={0021-2172}, mr={1045372}}
\bptok{imsref}%
\end{barticle}
\endbibitem

%b14 ###
\bibitem{Monroe72}
\begin{barticle}[mr]
\bauthor{\bsnm{Monroe},~\bfnm{Itrel}\binits{I.}}
(\byear{1972}).
\btitle{On embedding right continuous martingales in {B}rownian motion}.
\bjournal{Ann. Math. Statist.}
\bvolume{43}
\bpages{1293--1311}.
\bid{issn={0003-4851}, mr={0343354}}
\bptok{imsref}%
\end{barticle}
\endbibitem

%b15 ###
\bibitem{Monroe78}
\begin{barticle}[mr]
\bauthor{\bsnm{Monroe},~\bfnm{Itrel}\binits{I.}}
(\byear{1978}).
\btitle{Processes that can be embedded in {B}rownian motion}.
\bjournal{Ann. Probab.}
\bvolume{6}
\bpages{42--56}.
\bid{mr={0455113}}
\bptok{imsref}%
\end{barticle}
\endbibitem

%b16 ###
\bibitem{Obloj04}
\begin{barticle}[mr]
\bauthor{\bsnm{Ob{\l}{\'o}j},~\bfnm{Jan}\binits{J.}}
(\byear{2004}).
\btitle{The {S}korokhod embedding problem and its offspring}.
\bjournal{Probab. Surv.}
\bvolume{1}
\bpages{321--390}.
\bid{doi={10.1214/154957804100000060}, issn={1549-5787}, mr={2068476}}
\bptok{imsref}%
\end{barticle}
\endbibitem

%b17 ###
\bibitem{Perkins86}
\begin{bincollection}[mr]
\bauthor{\bsnm{Perkins},~\bfnm{Edwin}\binits{E.}}
(\byear{1986}).
\btitle{The {C}ereteli--{D}avis solution to the {$H\sp 1$}-embedding problem and
  an optimal embedding in {B}rownian motion}.
In \bbooktitle{Seminar on Stochastic Processes, 1985 ({G}ainesville, {F}la.,
  1985)}.
\bseries{Progr. Probab. Statist.}
\bvolume{12}
\bpages{172--223}.
\bpublisher{Birkh\"auser}, \blocation{Boston, MA}.
\bid{mr={0896743}}
\bptok{imsref}%
\end{bincollection}
\endbibitem

%b18 ###
\bibitem{RevuzYor99}
\begin{bbook}[mr]
\bauthor{\bsnm{Revuz},~\bfnm{Daniel}\binits{D.}} \AND
  \bauthor{\bsnm{Yor},~\bfnm{Marc}\binits{M.}}
(\byear{1999}).
\btitle{Continuous Martingales and {B}rownian Motion},
\bedition{3rd} ed.
\bseries{Grundlehren der Mathematischen Wissenschaften [Fundamental Principles
  of Mathematical Sciences]}
\bvolume{293}.
\bpublisher{Springer}, \blocation{Berlin}.
\bid{mr={1725357}}
\bptok{imsref}%
\end{bbook}
\endbibitem

%b19 ###
\bibitem{Rogers89}
\begin{barticle}[mr]
\bauthor{\bsnm{Rogers},~\bfnm{L.~C.~G.}\binits{L.~C.~G.}}
(\byear{1989}).
\btitle{A guided tour through excursions}.
\bjournal{Bull. Lond. Math. Soc.}
\bvolume{21}
\bpages{305--341}.
\bid{doi={10.1112/blms/21.4.305}, issn={0024-6093}, mr={0998631}}
\bptok{imsref}%
\end{barticle}
\endbibitem

%b20 ###
\bibitem{Rogers93}
\begin{barticle}[mr]
\bauthor{\bsnm{Rogers},~\bfnm{L.~C.~G.}\binits{L.~C.~G.}}
(\byear{1993}).
\btitle{The joint law of the maximum and terminal value of a martingale}.
\bjournal{Probab. Theory Related Fields}
\bvolume{95}
\bpages{451--466}.
\bid{doi={10.1007/BF01196729}, issn={0178-8051}, mr={1217446}}
\bptok{imsref}%
\end{barticle}
\endbibitem

%b21 ###
\bibitem{Skorokhod65}
\begin{bbook}[mr]
\bauthor{\bsnm{Skorokhod},~\bfnm{A.~V.}\binits{A.~V.}}
(\byear{1965}).
\btitle{Studies in the Theory of Random Processes}.
\bpublisher{Addison-Wesley}, \blocation{Reading, MA}.
\bid{mr={0185620}}
\bptok{imsref}%
\end{bbook}
\endbibitem

%b22 ###
\bibitem{Strassen67}
\begin{bincollection}[mr]
\bauthor{\bsnm{Strassen},~\bfnm{Volker}\binits{V.}}
(\byear{1967}).
\btitle{Almost sure behavior of sums of independent random variables and
  martingales}.
In \bbooktitle{Proc. {F}ifth {B}erkeley {S}ympos. {M}ath. {S}tatist. and
  {P}robability ({B}erkeley, {C}alif., 1965/66)}
\bpages{Vol. II: Contributions to Probability Theory, Part 1, pp. 315--343}.
\bpublisher{Univ. California Press}, \blocation{Berkeley, CA}.
\bid{mr={0214118}}
\bptok{imsref}%
\end{bincollection}
\endbibitem

%b23 ###
\bibitem{Vallois94}
\begin{barticle}[mr]
\bauthor{\bsnm{Vallois},~\bfnm{P.}\binits{P.}}
(\byear{1994}).
\btitle{Sur la loi du maximum et du temps local d'une martingale continue
  uniformement int\'egrable}.
\bjournal{Proc. Lond. Math. Soc. (3)}
\bvolume{69}
\bpages{399--427}.
\bid{doi={10.1112/plms/s3-69.2.399}, issn={0024-6115}, mr={1281971}}
\bptok{imsref}%
\end{barticle}
\endbibitem

\end{thebibliography}
\end{document}